# Polylinear Mapping of Free Algebra

Aleks Kleyn

Abstract. In this paper I consider the structure of the polylinear mapping of the free algebra over the commutative ring.

Contents

## 1. Algebra over Ring

**Definition 1.1.** Let $D$ be commutative ring. Let $A$ be module over ring $D$.[1] For given bilinear mapping

$$f : A \times A \to A$$

we define product in $A$

$$ab = f \circ (a, b) \tag{1.1}$$

$A$ is a **algebra over ring** $D$ if $A$ is $D$-module and we defined product (1.1) in $A$. If $A$ is free $D$-module, then $A$ is called **free algebra over ring** $D$. □

Aleks_Kleyn@MailAPS.org.
http://sites.google.com/site/AleksKleyn/.
http://arxiv.org/a/kleyn_a_1.
http://AleksKleyn.blogspot.com/.

[1] This subsection is written on the base of the section [4]-2.2.

Let $\overline{\overline{e}}$ be the basis of free algebra $A$ over ring $D$. If algebra $A$ has unit, then we assume that $\overline{e}_0$ is the unit of algebra $A$.

**Theorem 1.2.** *Let $\overline{\overline{e}}$ be the basis of free algebra $A$ over ring $D$. Let*

$$a = a^i e_i \quad b = b^i e_i \quad a, b \in A$$

*We can get the product of $a$, $b$ according to rule*

$$(ab)^k = C^k_{ij} a^i b^j \tag{1.2}$$

*where* $C^k_{ij}$ *are* **structural constants of algebra** $A$ **over ring** $D$. *The product of basis vectors in the algebra $A$ is defined according to rule*

$$\overline{e}_i \overline{e}_j = C^k_{ij} \overline{e}_k \tag{1.3}$$

*Proof.* The equation (1.3) is corollary of the statement that $\overline{\overline{e}}$ is the basis of the algebra $A$. Since the product in the algebra is a bilinear mapping, than we can write the product of $a$ and $b$ as

$$ab = a^i b^j e_i e_j \tag{1.4}$$

From equations (1.3), (1.4), it follows that

$$ab = a^i b^j C^k_{ij} \overline{e}_k \tag{1.5}$$

Since $\overline{\overline{e}}$ is a basis of the algebra $A$, than the equation (1.2) follows from the equation (1.5). □

**Definition 1.3.** Let $A_1$ and $A_2$ be algebras over ring $D$. The linear mapping[2]

$$f : A_1 \to A_2$$

of the $D$-module $A_1$ into the $D$-module $A_2$ is called **linear mapping of $D$-algebra $A_1$ into $D$-algebra** $A_2$. Let us denote $\mathcal{L}(A_1; A_2)$ set of linear mappings of algebra $A_1$ into algebra $A_2$. □

**Definition 1.4.** Let $D$ be the commutative associative ring. Let $A_1$, ..., $A_n$ be $D$-algebras and $S$ be $D$-module. We call map

$$f : A_1 \times ... \times A_n \to S$$

**polylinear mapping of algebras** $A_1$, ..., $A_n$ into module $S$, if

$$f \circ (a_1, ..., a_i + b_i, ..., a_n) = f \circ (a_1, ..., a_i, ..., a_n) + f \circ (a_1, ..., b_i, ..., a_n)$$
$$f \circ (a_1, ..., pa_i, ..., a_n) = pf \circ (a_1, ..., a_i, ..., a_n)$$

$$1 \le i \le n \quad a_i, b_i \in A_i \quad p \in D$$

Let us denote $\mathcal{L}(A_1, ..., A_n; S)$ set of polylinear mappings of algebras $A_1$, ..., $A_n$ into module $S$. Let us denote $\mathcal{L}(A^n; S)$ set of $n$-linear mappings of algebra $A$ ($A_1 = ... = A_n = A$) into module $S$. □

**Theorem 1.5.** *Let $D$ be the commutative associative ring. Let $A_1$, ..., $A_n$ be $D$-algebras and $S$ be $D$-module. Let mappings*

$$f : A_1 \times ... \times A_n \to S$$
$$g : A_1 \times ... \times A_n \to S$$

[2]This subsection is written on the base of the section [4]-2.3.

*be polylinear mappings. Then mapping $f+g$ defined by equation*

$$(f+g)\circ(a_1,...,a_n)=f\circ(a_1,...,a_n)+g\circ(a_1,...,a_n)$$

*is polylinear.*

*Proof.* Statement of theorem follows from chains of equations

$$\begin{aligned}&(f+g)\circ(x_1,...,x_i+y_i,...,x_n)\\=&f\circ(x_1,...,x_i+y_i,...,x_n)+g\circ(x_1,...,x_i+y_i,...,x_n)\\=&f\circ(x_1,...,x_i,...,x_n)+f\circ(x_1,...,y_i,...,x_n)\\+&g\circ(x_1,...,x_i,...,x_n)+g\circ(x_1,...,y_i,...,x_n)\\=&(f+g)\circ(x_1,...,x_i,...,x_n)+(f+g)\circ(x_1,...,y_i,...,x_n)\end{aligned}$$

$$\begin{aligned}&(f+g)\circ(x_1,...,px_i,...,x_n)\\=&f\circ(x_1,...,px_i,...,x_n)+g\circ(x_1,...,px_i,...,x_n)\\=&pf\circ(x_1,...,x_i,...,x_n)+pg\circ(x_1,...,x_i,...,x_n)\\=&p(f\circ(x_1,...,x_i,...,x_n)+g\circ(x_1,...,x_i,...,x_n))\\=&p(f+g)\circ(x_1,...,x_i,...,x_n)\end{aligned}$$

□

**Corollary 1.6.** *Consider algebra $A_1$ and algebra $A_2$. Let mappings*

$$f:A_1\to A_2$$
$$g:A_1\to A_2$$

*be linear mappings. Then mapping $f+g$ defined by equation*

$$(f+g)\circ a=f\circ a+g\circ a$$

*is linear.* □

**Theorem 1.7.** *Let $D$ be the commutative associative ring. Let $A_1$, ..., $A_n$ be $D$-algebras and $S$ be $D$-module. Let mapping*

$$f:A_1\times...\times A_n\to S$$

*be polylinear mapping. Then mapping $pf$, $p\in D$, defined by equation*

$$(pf)\circ x=p\ f\circ x$$

*is polylinear. This holds*

$$p(qf)=(pq)f$$
$$(p+q)f=pf+qf$$

*Proof.* Statement of theorem follows from chains of equations

$$\begin{aligned}
(pf)\circ(x_1,...,x_i+y_i,...,x_n)=&p\ f\circ(x_1,...,x_i+y_i,...,x_n)\\
=&p\ (f\circ(x_1,...,x_i,...,x_n)+f\circ(x_1,...,y_i,...,x_n))\\
=&p\ f\circ(x_1,...,x_i,...,x_n)+p\ f\circ(x_1,...,y_i,...,x_n)\\
=&(pf)\circ(x_1,...,x_i,...,x_n)+(pf)\circ(x_1,...,y_i,...,x_n)\\
(pf)\circ(x_1,...,qx_i,...,x_n)=&p\ f\circ(x_1,...,qx_i,...,x_n)=pq\ f\circ(x_1,...,x_i,...,x_n)\\
=&qp\ f\circ(x_1,...,x_n)=q\ (pf)\circ(x_1,...,x_n)\\
(p(qf))\circ(x_1,...,x_n)=&p\ (qf)\circ(x_1,...,x_n)=p\ (q\ f\circ(x_1,...,x_n))\\
=&(pq)\ f\circ(x_1,...,x_n)=((pq)f)\circ(x_1,...,x_n)\\
((p+q)f)\circ(x_1,...,x_n)=&(p+q)\ f\circ(x_1,...,x_n)\\
=&p\ f\circ(x_1,...,x_n)+q\ f\circ(x_1,...,x_n)\\
=&(pf)\circ(x_1,...,x_n)+(qf)\circ(x_1,...,x_n)
\end{aligned}$$

□

**Corollary 1.8.** *Consider algebra $A_1$ and algebra $A_2$. Let mapping*

$$f:A_1\to A_2$$

*be linear mapping. Then mapping $pf$, $p\in D$, defined by equation*

$$(pf)\circ x=p\ f\circ x$$

*is linear. This holds*

$$p(qf)=(pq)f$$
$$(p+q)f=pf+qf$$

□

**Theorem 1.9.** *Let $D$ be the commutative associative ring. Let $A_1$, ..., $A_n$ be $D$-algebras and $S$ be $D$-module. The set $\mathcal{L}(A_1,...,A_n;S)$ is a $D$-module.*

*Proof.* The theorem 1.5 determines the sum of polylinear mappings into $D$-module $S$. Let $f, g, h\in\mathcal{L}(A_1,...,A_2;S)$. For any $a=(a_1,...,a_n)$, $a_1\in A_1$, ..., $a_n\in A_n$,

$$\begin{aligned}
(f+g)\circ a=&f\circ a+g\circ a=g\circ a+f\circ a\\
=&(g+f)\circ a\\
((f+g)+h)\circ a=&(f+g)\circ a+h\circ a=(f\circ a+g\circ a)+h\circ a\\
=&f\circ a+(g\circ a+h\circ a)=f\circ a+(g+h)\circ a\\
=&(f+(g+h))\circ a
\end{aligned}$$

Therefore, sum of polylinear mappings is commutative and associative.

The mapping $z$ defined by equation

$$z\circ a=0$$

is zero of addition, because

$$(z+f)\circ a=z\circ a+f\circ a=0+f\circ a=f\circ a$$

For a given mapping $f$ a mapping $g$ defined by equation

$$g\circ a=-f\circ a$$

satisfies to equation

$$f + g = z$$

because

$$(f + g) \circ a = f \circ a + g \circ a = f \circ a - f \circ a = 0$$

Therefore, the set $\mathcal{L}(A_1; A_2)$ is an Abelian group.

From the theorem 1.7, it follows that the representation of the ring $D$ in the Abelian group $\mathcal{L}(A_1, ..., A_n; S)$ is defined. Since the ring $D$ has unit, than, according to the theorem [4]-2.1.1, specified representation is effective. □

**Corollary 1.10.** *Let $D$ be commutative ring with unit. Consider $D$-algebra $A_1$ and $D$-algebra $A_2$. The set $\mathcal{L}(A_1; A_2)$ is an $D$-module.* □

**Theorem 1.11.** *Let $A$ be algebra over commutative ring $D$. Module $\mathcal{L}(A; A)$ equiped by product*

$$\circ : (g, f) \in \mathcal{L}(A; A) \times \mathcal{L}(A; A) \to g \circ f \in \mathcal{L}(A; A) \tag{1.6}$$

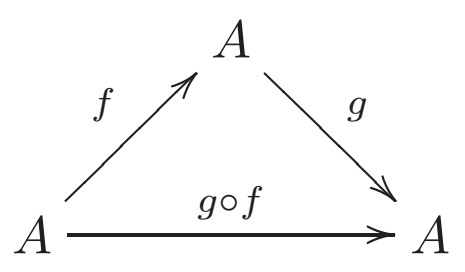

*is algebra over $D$.*

*Proof.* See the proof of the theorem [4]-2.4.5. □

## 2. Tensor Product of Algebras

**Definition 2.1.** Let $A_1$, ..., $A_n$ be free algebras over commutative ring $D$.[3] Let us consider category $\mathcal{A}$ whose objects are polylinear over commutative ring $D$ mappings

$$f : A_1 \times ... \times A_n \longrightarrow S_1 \qquad g : A_1 \times ... \times A_n \longrightarrow S_2$$

where $S_1$, $S_2$ are modules over ring $D$, We define morphism $f \to g$ to be linear over commutative ring $D$ mapping $h : S_1 \to S_2$ making diagram commutative. Universal object $A_1 \otimes ... \otimes A_n$ of category $\mathcal{A}$ is called **tensor product of algebras** $A_1$, ..., $A_n$. □

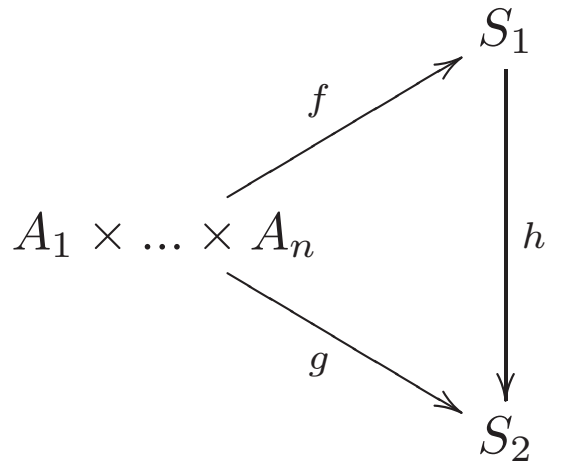

**Definition 2.2.** Tensor product

$$A^{\otimes n} = A_1 \otimes ... \otimes A_n \qquad A_1 = ... = A_n = A$$

is called **tensor power of algebra** $A$. □

[3] I give definition of tensor product of algebras following to definition in [1], p. 601 - 603. This subsection is written on the base of the section [4]-2.5.

**Theorem 2.3.** *There exists tensor product of algebras.*

*Proof.* Let $M$ be module over ring $D$ generated by product $A_1\times...\times A_n$ of algebras $A_1$, ..., $A_n$. Injection

$$i:A_1\times...\times A_n\longrightarrow M$$

is defined according to rule

$$i\circ(d_1,...,d_n)=(d_1,...,d_n) \tag{2.1}$$

Let $N\subset M$ be submodule generated by elements of the following type

$$(d_1,...,d_i+c_i,...,d_n)-(d_1,...,d_i,...,d_n)-(d_1,...,c_i,...,d_n) \tag{2.2}$$

$$(d_1,...,ad_i,...,d_n)-a(d_1,...,d_i,...,d_n) \tag{2.3}$$

where $d_i\in A_i$, $c_i\in A_i$, $a\in D$. Let

$$j:M\to M/N$$

be canonical map on factor module. Consider commutative diagram

(2.4)

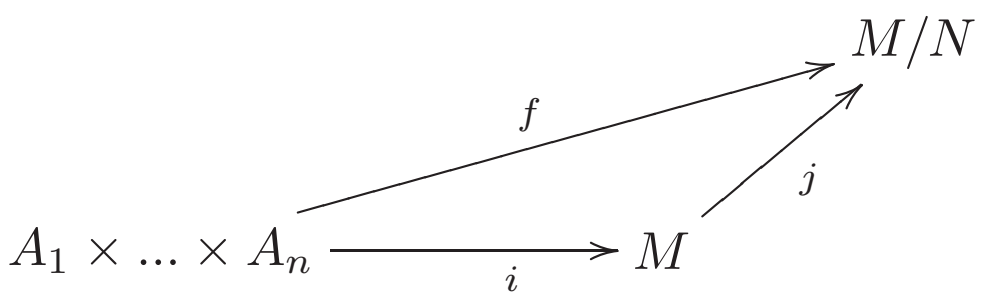

Since elements (2.2) and (2.3) belong to kernel of linear map $j$, then, from equation (2.1), it follows

$$f\circ(d_1,...,d_i+c_i,...,d_n)=f\circ(d_1,...,d_i,...,d_n)+f\circ(d_1,...,c_i,...,d_n) \tag{2.5}$$

$$f\circ(d_1,...,ad_i,...,d_n)=a\ f\circ(d_1,...,d_i,...,d_n) \tag{2.6}$$

From equations (2.5) and (2.6) it follows that map $f$ is polylinear over ring $D$. Since $M$ is module with basis $A_1\times...\times A_n$, than, according to theorem [1]-4.1 on p. 135, for any module $V$ and any polylinear over $D$ map

$$g:A_1\times...\times A_n\longrightarrow V$$

there exists a unique homomorphism $k:M\to V$, for which following diagram is commutative

(2.7)

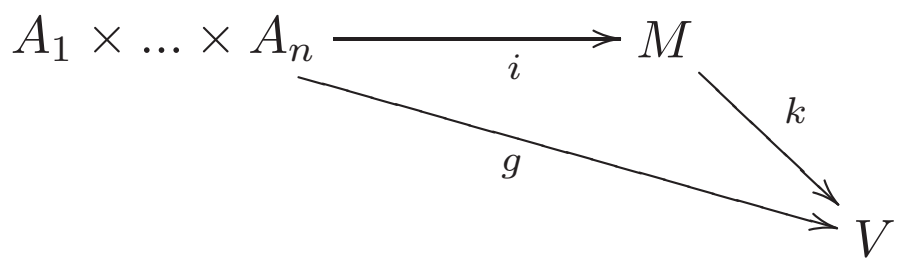

Since $g$ is polylinear over $D$, than $\ker k\subseteq N$. According to statement on p. [1]-119, map $j$ is universal in the category of homomorphisms of vector space $M$ whose kernel contains $N$. Therefore, we have homomorphism

$$h:M/N\to V$$

which makes the following diagram commutative

(2.8)

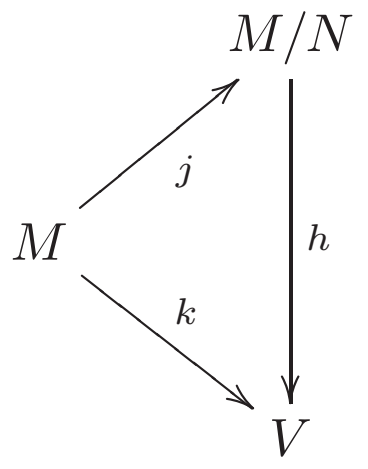

We join diagrams (2.4), (2.7), (2.8), and get commutative diagram

(2.9)

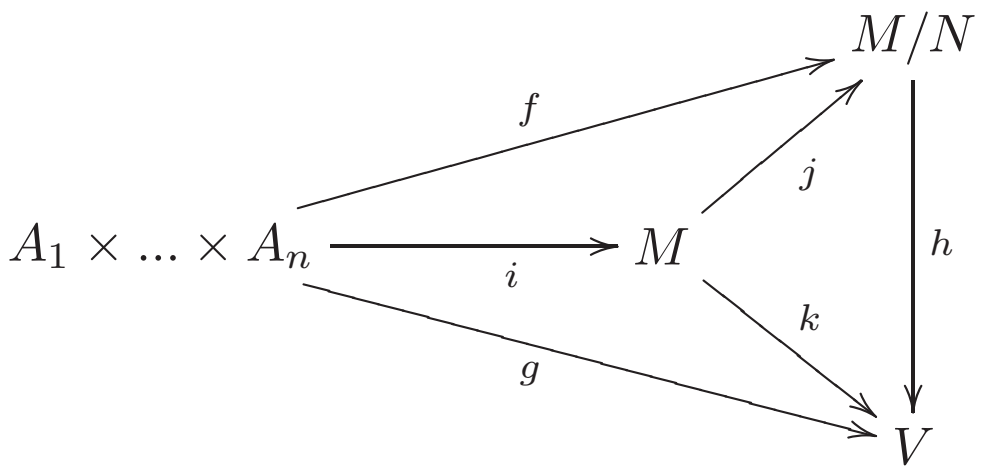

Since Im$f$ generates $M/N$, than map $h$ is uniquely determined. □

According to proof of theorem 2.3

$$A_1 \otimes ... \otimes A_n = M/N$$

If $d_i \in A_i$, we write

(2.10) $$j \circ (d_1, ..., d_n) = d_1 \otimes ... \otimes d_n$$

**Theorem 2.4.** *Let $A_1$, ..., $A_n$ be algebras over commutative ring $D$. Let*

$$f : A_1 \times ... \times A_n \to A_1 \otimes ... \otimes A_n$$

*be polylinear mapping defined by equation*

(2.11) $$f \circ (d_1, ..., d_n) = d_1 \otimes ... \otimes d_n$$

*Let*

$$g : A_1 \times ... \times A_n \to V$$

*be polylinear mapping into $D$-module $V$. There exists an $D$-linear mapping*

$$h : A_1 \otimes ... \otimes A_n \to V$$

*such that the diagram*

(2.12)

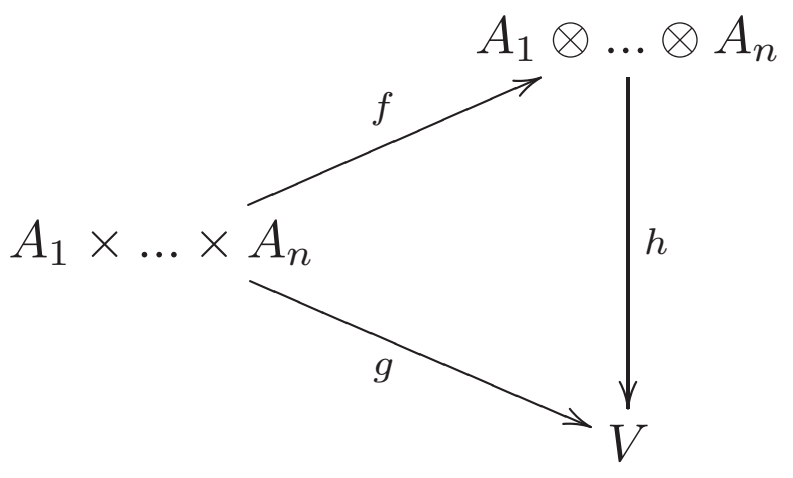

*is commutative.*

*Proof.* Equation (2.11) follows from equations (2.1) and (2.10). An existence of the mapping $h$ follows from the definition 2.1 and constructions made in the proof of the theorem 2.3. □

We can write equations (2.5) and (2.6) as

$$\begin{aligned}&a_1\otimes...\otimes(a_i+b_i)\otimes...\otimes a_n\\=&a_1\otimes...\otimes a_i\otimes...\otimes a_n+a_1\otimes...\otimes b_i\otimes...\otimes a_n\end{aligned}\tag{2.13}$$

$$a_1\otimes...\otimes(ca_i)\otimes...\otimes a_n=c(a_1\otimes...\otimes a_i\otimes...\otimes a_n)\tag{2.14}$$

$$a_i\in A_i\quad b_i\in A_i\quad c\in D$$

**Theorem 2.5.** *Let $A$ be algebra over commutative ring $D$. There exists a linear mapping*

$$h:a\otimes b\in A\otimes A\rightarrow ab\in A$$

*Proof.* The theorem is corollary of the theorem 2.4 and the definition 1.1. □

**Theorem 2.6.** *Tensor product $A_1\otimes...\otimes A_n$ of free finite dimensional algebras $A_1$, ..., $A_n$ over the commutative ring $D$ is free finite dimensional algebra.*

*Let $\overline{\overline{e}}_i$ be the basis of algebra $A_i$ over ring $D$. We can represent any tensor $a\in A_1\otimes...\otimes A_n$ in the following form*

$$a=a^{i_1...i_n}\overline{e}_{1\cdot i_1}\otimes...\otimes\overline{e}_{n\cdot i_n}\tag{2.15}$$

*Expression $a^{i_1...i_n}$ is called* **standard component of tensor**.

*Proof.* See the proof of the theorem [4]-2.5.6. □

## 3. Linear Mapping into Associative Algebra

**Theorem 3.1.** *Consider $D$-algebras $A_1$ and $A_2$. For given mapping $f\in\mathcal{L}(A_1;A_2)$, there exists linear mapping*

$$h:A_2\otimes A_2\rightarrow\mathcal{L}(A_1;A_2)$$

*defined by the equation*

$$(a\otimes b)\circ f=afb\tag{3.1}$$

*Proof.* See the proof of the theorems [4]-2.6.1 and [4]-2.6.2. □

**Theorem 3.2.** *Consider $D$-algebras $A_1$ and $A_2$. Let us define product in algebra $A_2\otimes A_2$ according to rule*

$$(c\otimes d)\circ(a\otimes b)=(ca)\otimes(bd)\tag{3.2}$$

*A linear mapping*

$$h:A_2\otimes A_2\rightarrow{}^*\mathcal{L}(A_1;A_2)\tag{3.3}$$

*defined by the equation*

$$(a\otimes b)\circ f=afb\quad a,b\in A_2\quad f\in\mathcal{L}(A_1;A_2)\tag{3.4}$$

*is representation*[4] *of algebra $A_2\otimes A_2$ in module $\mathcal{L}(A_1;A_2)$.*

*Proof.* See the proof of the theorem [4]-2.6.3. □

[4]See the definition of representation of $\Omega$-algebra in the definition [3]-2.1.4.

**Theorem 3.3.** *Consider $D$-algebra $A$. Let us define product in algebra $A\otimes A$ according to rule (3.2). A representation of algebra $A\otimes A$*

$$h:A\otimes A\rightarrow {}^*\mathcal{L}(A;A) \tag{3.5}$$

*in module $\mathcal{L}(A;A)$ defined by the equation*

$$(a\otimes b)\circ f=afb\quad a,b\in A\quad f\in\mathcal{L}(A;A) \tag{3.6}$$

*allows us to identify tensor $d\in A\otimes A$ and mapping $d\circ\delta\in\mathcal{L}(A;A)$ where $\delta\in\mathcal{L}(A;A)$ is identity mapping.*

*Proof.* See the proof of the theorem [4]-2.6.4. □

From the theorem 3.3, it follows that we can consider the mapping (3.4) as the product of mappings $a\otimes b$ and $f$.

**Theorem 3.4.** *Consider $D$-algebra $A_1$ and associative $D$-algebra $A_2$. Consider the representation of algebra $A_2\otimes A_2$ in the module $\mathcal{L}(A_1;A_2)$. The mapping*

$$h:A_1\rightarrow A_2$$

*generated by the mapping*

$$f:A_1\rightarrow A_2$$

*has form*

$$h=(a_{s\cdot 0}\otimes a_{s\cdot 1})\circ f=a_{s\cdot 0}fa_{s\cdot 1} \tag{3.7}$$

*Proof.* See the proof of the theorem [4]-2.6.6 □

**Theorem 3.5.** *Let $A$ be free finite dimensional associative algebra over the field $D$. The representation of algebra $A\otimes A$ in algebra $\mathcal{L}(A;A)$ has finite* **basis** *$\overline{\overline{I}}$.*

(1) *The linear mapping $f\in\mathcal{L}(A;A)$ has form*

$$f=(a_{k\cdot s_k\cdot 0}\otimes a_{k\cdot s_k\cdot 1})\circ I_k=\sum_{k}a_{k\cdot s_k\cdot 0}I_ka_{k\cdot s_k\cdot 1} \tag{3.8}$$

(2) *Its standard representation has form*

$$f=a^{k\cdot ij}(\overline{e}_i\otimes\overline{e}_j)\circ I_k=a^{k\cdot ij}\overline{e}_iI_k\overline{e}_j \tag{3.9}$$

*Proof.* See the proof of the theorem [4]-2.7.5 □

## 4. Linear Mapping into Nonassociative Algebra

Since the product is nonassociative, we may assume that action of $a,b\in A$ over the mapping $f$ may have form either $a(fb)$, or $(af)b$.

**Theorem 4.1.** *Let $\overline{\overline{e}}_1$ be basis of the free finite dimensional $D$-algebra $A_1$. Let $\overline{\overline{e}}_2$ be basis of the free finite dimensional nonassociative $D$-algebra $A_2$. Let $C_{2\cdot}{}^{p}_{kl}$ be structural constants of algebra $A_2$. Let the mapping*

$$g=a\circ f \tag{4.1}$$

*generated by the mapping $f\in(A_1;A_2)$ through the tensor $a\in A_2\otimes A_2$, has the standard representation*

$$g=a^{ij}(\overline{e}_i\otimes\overline{e}_j)\circ f=a^{ij}(\overline{e}_if)\overline{e}_j \tag{4.2}$$

*Coordinates of the mapping (4.1) and its standard components are connected by the equation*

$$g_l^k = f_l^m g^{ij} C_{2\cdot}{}^p_{im} C_{2\cdot}{}^k_{pj} \tag{4.3}$$

*Proof.* See the proof of the theorem [4]-2.8.3 □

**Theorem 4.2.** *Let $A$ be free finite dimensional nonassociative algebra over the ring $D$. The representation of algebra $A\otimes A$ in algebra $\mathcal{L}(A;A)$ has finite basis $\overline{\overline{I}}$.*

(1) *The linear mapping $f\in\mathcal{L}(A;A)$ has form*

$$f = (a_{k\cdot s_k\cdot 0}\otimes a_{k\cdot s_k\cdot 1})\circ I_k = (a_{k\cdot s_k\cdot 0} I_k) a_{k\cdot s_k\cdot 1} \tag{4.4}$$

(2) *Its standard representation has form*

$$f = a^{k\cdot ij}(\overline{e}_i\otimes\overline{e}_j)\circ I_k = a^{k\cdot ij}(\overline{e}_i I_k)\overline{e}_j \tag{4.5}$$

*Proof.* See the proof of the theorem [4]-2.8.4 □

## 5. Polylinear Mapping into Associative Algebra

**Theorem 5.1.** *Let $A_1$, ..., $A_n$, $A$ be associative $D$-algebras. Let*

$$f_i\in\mathcal{L}(A_i;A)\quad i=1,...,n$$

$$a_j\in A\quad j=0,...,n$$

*For given transposition $\sigma$ of $n$ variables, the mapping*

$$\begin{aligned}
&((a_0,...,a_n)\circ\sigma(f_1,...,f_n))\circ(x_1,...,x_n)\\
=&\ (a_0\sigma(f_1)a_1...a_{n-1}\sigma(f_n)a_n)\circ(x_1,...,x_n)\\
=&\ a_0\sigma(f_1\circ x_1)a_1...a_{n-1}\sigma(f_n\circ x_n)a_n
\end{aligned}\tag{5.1}$$

*is $n$-linear mapping into algebra $A$.*

*Proof.* The statement of theorem follows from chains of equations

$$\begin{aligned}
&((a_0,...,a_n,\sigma)\circ(f_1,...,f_n))\circ(x_1,...,x_i+y_i,...,x_n)\\
=&a_0\sigma(f_1\circ x_1)a_1...\sigma(f_i\circ(x_i+y_i))...a_{n-1}\sigma(f_n\circ x_n)a_n\\
=&a_0\sigma(f_1\circ x_1)a_1...\sigma(f_i\circ x_i+f_i\circ y_i)...a_{n-1}\sigma(f_n\circ x_n)a_n\\
=&a_0\sigma(f_1\circ x_1)a_1...\sigma(f_i\circ x_i)...a_{n-1}\sigma(f_n\circ x_n)a_n\\
+&a_0\sigma(f_1\circ x_1)a_1...\sigma(f_i\circ y_i)...a_{n-1}\sigma(f_n\circ x_n)a_n\\
=&((a_0,...,a_n,\sigma)\circ(f_1,...,f_n))\circ(x_1,...,x_i,...,x_n)\\
+&((a_0,...,a_n,\sigma)\circ(f_1,...,f_n))\circ(x_1,...,y_i,...,x_n)
\end{aligned}$$

$$\begin{aligned}
&((a_0,...,a_n,\sigma)\circ(f_1,...,f_n))\circ(x_1,...,px_i,...,x_n)\\
=&a_0\sigma(f_1\circ x_1)a_1...\sigma(f_i\circ(px_i))...a_{n-1}\sigma(f_n\circ x_n)a_n\\
=&a_0\sigma(f_1\circ x_1)a_1...\sigma(p(f_i\circ x_i))...a_{n-1}\sigma(f_n\circ x_n)a_n\\
=&p(a_0\sigma(f_1\circ x_1)a_1...\sigma(f_i\circ x_i)...a_{n-1}\sigma(f_n\circ x_n)a_n)\\
=&p(((a_0,...,a_n,\sigma)\circ(f_1,...,f_n))\circ(x_1,...,x_i,...,x_n))
\end{aligned}$$

□

In the equation (5.1), as well as in other expressions of polylinear mapping, we have convention that mapping $f_i$ has variable $x_i$ as argument.

**Theorem 5.2.** *Let $A_1$, ..., $A_n$, $A$ be associative $D$-algebras. For given set of mappings*

$$f_i\in\mathcal{L}(A_i;A)\quad i=1,...,n$$

*the mapping*

$$h:A^{n+1}\to\mathcal{L}(A_1,...,A_n;A)$$

*defined by equation*

$$(a_0,...,a_n,\sigma)\circ(f_1,...,f_n)=a_0\sigma(f_1)a_1...a_{n-1}\sigma(f_n)a_n$$

*is $n+1$-linear mapping into $D$-module $\mathcal{L}(A_1,...,A_n;A)$.*

*Proof.* The statement of theorem follows from chains of equations

$$\begin{aligned}&((a_0,...,a_i+b_i,...a_n,\sigma)\circ(f_1,...,f_n))\circ(x_1,...,x_n)\\
=&a_0\sigma(f_1\circ x_1)a_1...(a_i+b_i)...a_{n-1}\sigma(f_n\circ x_n)a_n\\
=&a_0\sigma(f_1\circ x_1)a_1...a_i...a_{n-1}\sigma(f_n\circ x_n)a_n+a_0\sigma(f_1\circ x_1)a_1...b_i...a_{n-1}\sigma(f_n\circ x_n)a_n\\
=&((a_0,...,a_i,...,a_n,\sigma)\circ(f_1,...,f_n))\circ(x_1,...,x_n)\\
+&((a_0,...,b_i,...,a_n,\sigma)\circ(f_1,...,f_n))\circ(x_1,...,x_n)\\
=&((a_0,...,a_i,...,a_n,\sigma)\circ(f_1,...,f_n)+(a_0,...,b_i,...,a_n,\sigma)\circ(f_1,...,f_n))\circ(x_1,...,x_n)\end{aligned}$$

$$\begin{aligned}&((a_0,...,pa_i,...a_n,\sigma)\circ(f_1,...,f_n))\circ(x_1,...,x_n)\\
=&a_0\sigma(f_1\circ x_1)a_1...pa_i...a_{n-1}\sigma(f_n\circ x_n)a_n\\
=&p(a_0\sigma(f_1\circ x_1)a_1...a_i...a_{n-1}\sigma(f_n\circ x_n)a_n)\\
=&p(((a_0,...,a_i,...,a_n,\sigma)\circ(f_1,...,f_n))\circ(x_1,...,x_n))\\
=&(p((a_0,...,a_i,...,a_n,\sigma)\circ(f_1,...,f_n)))\circ(x_1,...,x_n)\end{aligned}$$

□

**Theorem 5.3.** *Let $A_1$, ..., $A_n$, $A$ be associative $D$-algebras. For given set of mappings*

$$f_i\in\mathcal{L}(A_i;A)\quad i=1,...,n$$

*there exists linear mapping*

$$h:A^{\otimes n+1}\times S_n\to\mathcal{L}(A_1,...,A_n;A)$$

*defined by the equation*

$$\begin{aligned}(a_0\otimes...\otimes a_n,\sigma)\circ(f_1,...,f_n)&=(a_0,...,a_n,\sigma)\circ(f_1,...,f_n)\\&=a_0\sigma(f_1)a_1...a_{n-1}\sigma(f_n)a_n\end{aligned}\tag{5.2}$$

*Proof.* The statement of the theorem is corollary of theorems 2.4, 5.2. □

**Theorem 5.4.** *Let $A_1$, ..., $A_n$, $A$ be associative $D$-algebras. For given tensor $a\in A^{\otimes n+1}$ and given transposition $\sigma\in S_n$ the mapping*

$$h:\prod_{i=1}^n\mathcal{L}(A_i;A)\to\mathcal{L}(A_1,...,A_n;A)$$

*defined by equation*

$$(a_0\otimes...\otimes a_n,\sigma)\circ(f_1,...,f_n)=a_0\sigma(f_1)a_1...a_{n-1}\sigma(f_n)a_n$$

*is $n$-linear mapping into $D$-module $\mathcal{L}(A_1,...,A_n;A)$.*

*Proof.* The statement of theorem follows from chains of equations

$$
\begin{aligned}
&((a_0\otimes...\otimes a_n,\sigma)\circ(f_1,...,f_i+g_i,...,f_n))\circ(x_1,...,x_n)\\
=&(a_0\sigma(f_1)a_1...\sigma(f_i+g_i)...a_{n-1}\sigma(f_n)a_n)\circ(x_1,...,x_n)\\
=&a_0\sigma(f_1\circ x_1)a_1...\sigma((f_i+g_i)\circ x_i)...a_{n-1}\sigma(f_n\circ x_n)a_n\\
=&a_0\sigma(f_1\circ x_1)a_1...\sigma(f_i\circ x_i+g_i\circ x_i)...a_{n-1}\sigma(f_n\circ x_n)a_n\\
=&a_0\sigma(f_1\circ x_1)a_1...\sigma(f_i\circ x_i)...a_{n-1}\sigma(f_n\circ x_n)a_n\\
+&a_0\sigma(f_1\circ x_1)a_1...\sigma(g_i\circ x_i)...a_{n-1}\sigma(f_n\circ x_n)a_n\\
=&(a_0\sigma(f_1)a_1...\sigma(f_i)...a_{n-1}\sigma(f_n)a_n)\circ(x_1,...,x_n)\\
+&(a_0\sigma(f_1)a_1...\sigma(g_i)...a_{n-1}\sigma(f_n)a_n)\circ(x_1,...,x_n)\\
=&((a_0\otimes...\otimes a_n,\sigma)\circ(f_1,...,f_i,...,f_n))\circ(x_1,...,x_n)\\
+&((a_0\otimes...\otimes a_n,\sigma)\circ(f_1,...,g_i,...,f_n))\circ(x_1,...,x_n)\\
=&((a_0\otimes...\otimes a_n,\sigma)\circ(f_1,...,f_i,...,f_n)\\
+&(a_0\otimes...\otimes a_n,\sigma)\circ(f_1,...,g_i,...,f_n))\circ(x_1,...,x_n)
\end{aligned}
$$

$$
\begin{aligned}
&((a_0\otimes...\otimes a_n,\sigma)\circ(f_1,...,pf_i,...,f_n))\circ(x_1,...,x_n)\\
=&(a_0\sigma(f_1)a_1,...\sigma(pf_i)...a_{n-1}\sigma(f_n)a_n)\circ(x_1,...,x_n)\\
=&a_0\sigma(f_1\circ x_1)a_1,...\sigma((pf_i)\circ x_i)...a_{n-1}\sigma(f_n\circ x_n)a_n\\
=&a_0\sigma(f_1\circ x_1)a_1,...\sigma(p(f_i\circ x_i))...a_{n-1}\sigma(f_n\circ x_n)a_n\\
=&p(a_0\sigma(f_1\circ x_1)a_1,...\sigma(f_i\circ x_i)...a_{n-1}\sigma(f_n\circ x_n)a_n)\\
=&p(((a_0\otimes...\otimes a_n,\sigma)\circ(f_1,...,f_i,...,f_n))\circ(x_1,...,x_n))\\
=&(p((a_0\otimes...\otimes a_n,\sigma)\circ(f_1,...,f_i,...,f_n)))\circ(x_1,...,x_n)
\end{aligned}
$$

□

**Theorem 5.5.** *Let $A_1$, ..., $A_n$, $A$ be associative $D$-algebras. For given tensor $a\in A^{\otimes n+1}$ and given transposition $\sigma\in S_n$ there exists linear mapping*

$$h:\mathcal{L}(A_1;A)\otimes...\otimes\mathcal{L}(A_n;A)\rightarrow\mathcal{L}(A_1,...,A_n;A)$$

*defined by the equation*

$$(a_0\otimes...\otimes a_n,\sigma)\circ(f_1\otimes...\otimes f_n)=(a_0\otimes...\otimes a_n,\sigma)\circ(f_1,...,f_n)\tag{5.3}$$

*Proof.* The statement of the theorem is corollary of theorems 2.4, 5.4. □

**Theorem 5.6.** *Let $A$ be associative $D$-algebra. Polylinear map*

$$f:A^n\rightarrow A, a=f\circ(a_1,...,a_n)\tag{5.4}$$

*has form*

$$a=f^n_{s\cdot 0}\ \sigma_s(I_{1\cdot s}\circ a_1)\ f^n_{s\cdot 1}\ ...\ \sigma_s(I_{n\cdot s}\circ a_n)\ f^n_{s\cdot n}\tag{5.5}$$

*where $\sigma_s$ is a transposition of set of variables $\{a_1,...,a_n\}$*

$$\sigma_s=\begin{pmatrix}a_1&...&a_n\\ \sigma_s(a_1)&...&\sigma_s(a_n)\end{pmatrix}$$

*Proof.* We prove statement by induction on $n$.

When $n = 1$ the statement of theorem is corollary of the statement (1) of the theorem 3.5. In such case we may identify[5]

$$f^1_{s\cdot p} = f_{s\cdot p} \quad p = 0, 1$$

Let statement of theorem be true for $n = k - 1$. Then it is possible to represent mapping (5.4) as

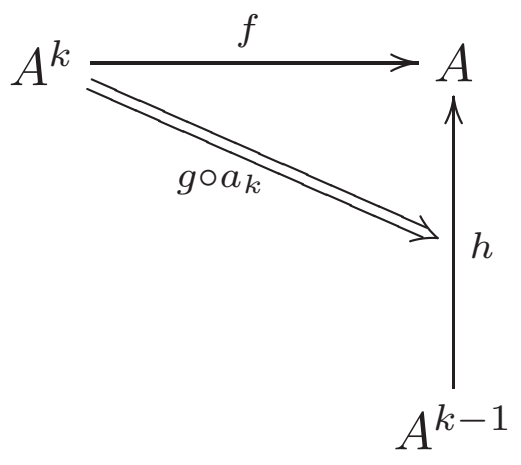

$$a = f \circ (a_1, ..., a_k) = (g \circ a_k) \circ (a_1, ..., a_{k-1})$$

According to statement of induction polylinear mapping $h$ has form

$$a = h^{k-1}_{t\cdot 0}\ \sigma_t(I_{1\cdot t} \circ a_1)\ h^{k-1}_{t\cdot 1}\ ...\ \sigma_t(I_{k-1\cdot t} \circ a_{k-1})\ h^{k-1}_{t\cdot k-1}$$

According to construction $h = g \circ a_k$. Therefore, expressions $h_{t\cdot p}$ are functions of $a_k$. Since $g \circ a_k$ is linear mapping of $a_k$, then only one expression $h_{t\cdot p}$ is linear mapping of $a_k$, and rest expressions $h_{t\cdot q}$ do not depend on $a_k$.

Without loss of generality, assume $p = 0$. According to the equation (3.7) for given $t$

$$h^{k-1}_{t\cdot 0} = g_{tr\cdot 0}\ I_{k\cdot r} \circ a_k\ g_{tr\cdot 1}$$

Assume $s = tr$. Let us define transposition $\sigma_s$ according to rule

$$\sigma_s = \sigma_{tr} = \begin{pmatrix} a_k & a_1 & ... & a_{k-1} \\ a_k & \sigma_t(a_1) & ... & \sigma_t(a_{k-1}) \end{pmatrix}$$

Suppose

$$f^k_{tr\cdot q+1} = h^{k-1}_{t\cdot q} \quad q = 1, ..., k-1$$
$$f^k_{tr\cdot q} = g_{tr\cdot q} \quad q = 0, 1$$

We proved step of induction. □

**Definition 5.7.** Expression $f^n_{s\cdot p}$ in equation (5.5) is called **component of polylinear map** $f$. □

**Theorem 5.8.** *Consider* $D$*-algebras* $A$*. A linear mapping*

$$h : A^{\otimes n+1} \times S_n \to {}^*\mathcal{L}(A^n; A)$$

*defined by the equation*

$$(a_0 \otimes ... \otimes a_n, \sigma) \circ (f_1 \otimes ... \otimes f_n) = a_0\sigma(f_1)a_1...a_{n-1}\sigma(f_n)a_n \tag{5.6}$$
$$a_0, ..., a_n \in A \quad \sigma \in S_n \quad f_1, ..., f_n \in \mathcal{L}(A_n; A)$$

[5]In representation (5.5) we will use following rules.

- If range of any index is set consisting of one element, then we will omit corresponding index.
- If $n = 1$, then $\sigma_s$ is identical transformation. We will not show such transformation in the expression.

*is representation*[6] *of algebra* $A^{\otimes n+1}\times S^n$ *in* $D$*-module* $\mathcal{L}(A^n;A)$*.*

*Proof.* According to the theorems 3.5, 5.6, we can represent $n$-linear mapping as sum of terms (5.1), where $f_i$, $i=1$, ..., $n$, are generators of representation (3.3). Let us write the term $s$ of the expression (5.5) as

$$b_1\sigma(I_{1\cdot s}\circ x_1)c_1b_2...c_{n-1}b_n\sigma(I_{n\cdot s}\circ x_n)c_n \tag{5.7}$$

where

$$b_1=f^n_{s\cdot 0}\quad b_2=...=b_n=e\quad c_1=f^n_{s\cdot 1}\quad ...\quad c_n=f^n_{s\cdot n}$$

Let us assume

$$f_i=\sigma^{-1}(b_i)I_{i\cdot s}\sigma^{-1}(c_i)\quad i=1,...,n$$

in equation (5.7). Therefore, according to theorem 5.3, mapping (5.6) is transformation of module $\mathcal{L}(A^n;A)$. For a given tensor $c\in A^{\otimes n+1}$ and given transposition $\sigma\in S_n$, a transformation $h(c,\sigma)$ is a linear transformation of module $\mathcal{L}(A^n;A)$ according to the theorem 5.5. According to the theorem 5.3, mapping (5.6) is linear mapping. According to the definition [3]-2.1.4 mapping (5.6) is a representation of the algebra $A^{\otimes n+1}\times S^n$ in the module $\mathcal{L}(A^n;A)$. □

**Theorem 5.9.** *Consider* $D$*-algebra* $A$*. A representation*

$$h:A^{\otimes n+1}\times S_n\rightarrow {}^*\mathcal{L}(A^n;A)$$

*of algebra* $A^{\otimes n+1}$ *in module* $\mathcal{L}(A^n;A)$ *defined by the equation* (5.6) *allows us to identify tensor* $d\in A^{\otimes n+1}$ *and transposition* $\sigma\in S^n$ *with mapping*

$$(d,\sigma)\circ(f_1,...,f_n)\quad f_i=\delta\in\mathcal{L}(A;A) \tag{5.8}$$

*where* $\delta\in\mathcal{L}(A;A)$ *is identity mapping.*

*Proof.* If we assume $f_i=\delta$, $d=a_0\otimes...\otimes a_n$ in the equation (5.2), then the equation (5.2) gets form

$$\begin{aligned}((a_0\otimes...\otimes a_n,\sigma)\circ(\delta,...,\delta))\circ(x_1,...,x_n)&=a_0\ (\delta\circ x_1)\ ...\ (\delta\circ x_n)\ a_n\\&=a_0\ x_1\ ...\ x_n\ a_n\end{aligned} \tag{5.9}$$

If we assume

$$\begin{aligned}&((a_0\otimes...\otimes a_n,\sigma)\circ(\delta,...,\delta))\circ(x_1,...,x_n)\\=&(a_0\otimes...\otimes a_n,\sigma)\circ(\delta\circ x_1,...,\delta\circ x_n)\\=&(a_0\otimes...\otimes a_n,\sigma)\circ(x_1,...,x_n)\end{aligned} \tag{5.10}$$

then comparison of equations (5.9) and (5.10) gives a basis to identify the action of the tensor $d=a_0\otimes...\otimes a_n$ and transposition $\sigma\in S^n$ with mapping (5.8). □

Instead of notation $(a_0\otimes...a_n,\sigma)$, we also use notation

$$a_0\otimes_{\sigma(1)}...\otimes_{\sigma(n)}a_n$$

when we want to show order of arguments in expression. For instance, the following expressions are equivalent

$$(a_0\otimes a_1\otimes a_2\otimes a_3,(2,1,3))\circ(x_1,x_2,x_3)=a_0x_2a_1x_1a_2x_3a_3$$

$$(a_0\otimes_2 a_1\otimes_1 a_2\otimes_3 a_3)\circ(x_1,x_2,x_3)=a_0x_2a_1x_1a_2x_3a_3$$

[6]See the definition of representation of $\Omega$-algebra in the definition [3]-2.1.4.

## 6. Polylinear Map into Free Finite Dimensional Associative Algebra

**Theorem 6.1.** *Let $A$ be free finite dimensional associative algebra over the ring $D$. Let $\overline{\overline{I}}$ be basis of algebra $\mathcal{L}(A;A)$. Let $\overline{\overline{e}}$ be the basis of the algebra $A$ over the ring $D$.* **Standard representation of polylinear mapping into associative algebra** *has form*

$$f\circ(a_1,...,a_n)=f^{i_0...i_n}_{t\cdot k_1...k_n}\ \overline{e}_{i_0}\ \sigma_t(I_{k_1}\circ a_1)\ \overline{e}_{i_1}\ ...\ \sigma_t(I_{k_n}\circ a_n)\ \overline{e}_{i_n} \tag{6.1}$$

*Index $t$ enumerates every possible transpositions $\sigma_t$ of the set of variables $\{a_1,...,a_n\}$. Expression $f^{i_0...i_n}_{t\cdot k_1...k_n}$ in equation (6.1) is called* **standard component of polylinear mapping** *$f$.*

*Proof.* We change index $s$ in the equation (5.5) so as to group the terms with the same set of generators $I_k$. Expression (5.5) gets form

$$a=f^n_{k_1...k_n\cdot s\cdot 0}\ \sigma_s(I_{k_1\cdot s}\circ a_1)\ f^n_{k_1...k_n\cdot s\cdot 1}\ ...\ \sigma_s(I_{k_n\cdot s}\circ a_n)\ f^n_{k_1...k_n\cdot s\cdot n} \tag{6.2}$$

We assume that the index $s$ takes values depending on $k_1$, ..., $k_n$. Components of polylinear map $f$ have expansion

$$f^n_{k_1...k_n\cdot s\cdot p}=\overline{e}_i f^{ni}_{k_1...k_n\cdot s\cdot p} \tag{6.3}$$

relative to basis $\overline{\overline{e}}$. If we substitute (6.3) into (5.5), we get

$$a=f^{nj_1}_{k_1...k_n\cdot s\cdot 0}\ \overline{e}_{j_1}\ \sigma_s(I_{k_1}\circ a_1)\ f^{nj_2}_{k_1...k_n\cdot s\cdot 1}\ \overline{e}_{j_2}\ ...\ \sigma_s(I_{k_n}\circ a_n)\ f^{nj_n}_{k_1...k_n\cdot s\cdot n}\ \overline{e}_{j_n} \tag{6.4}$$

Let us consider expression

$$f^{j_0...j_n}_{t\cdot k_1...k_n}=f^{nj_1}_{k_1...k_n\cdot s\cdot 0}\ ...f^{nj_n}_{k_1...k_n\cdot s\cdot n} \tag{6.5}$$

The right-hand side is supposed to be the sum of the terms with the index $s$, for which the transposition $\sigma_s$ is the same. Each such sum has a unique index $t$. If we substitute expression (6.5) into equation (6.4) we get equation (6.1). □

**Theorem 6.2.** *Let $A$ be free finite dimensional associative algebra over the ring $D$. Let $\overline{\overline{e}}$ be the basis of the algebra $A$ over the ring $D$. Polylinear map (5.4) can be represented as $D$-valued form of degree $n$ over ring $D$*[7]

$$f(a_1,...,a_n)=a_1^{i_1}...a_n^{i_n}f_{i_1...i_n} \tag{6.6}$$

*where*

$$\begin{aligned}a_j&=\overline{e}_i a^i_j\\ f_{i_1...i_n}&=f\circ(\overline{e}_{i_1},...,\overline{e}_{i_n})\end{aligned} \tag{6.7}$$

*and values $f_{i_1...i_n}$ are coordinates of $D$-valued covariant tensor.*

*Proof.* According to the definition 1.4, the equation (6.6) follows from the chain of equations

$$f\circ(a_1,...,a_n)=f\circ(\overline{e}_{i_1}a_1^{i_1},...,\overline{e}_{i_n}a_n^{i_n})=a_1^{i_1}...a_n^{i_n}f\circ(\overline{e}_{i_1},...,\overline{e}_{i_n})$$

Let $\overline{\overline{e}}'$ be another basis. Let

$$\overline{e}'_i=\overline{e}_j h^j_i \tag{6.8}$$

[7]We proved the theorem by analogy with theorem in [2], p. 107, 108

be transformation, mapping basis $\overline{\overline{e}}$ into basis $\overline{\overline{e}}'$. From equations (6.8) and (6.7) it follows

(6.9)
$$\begin{aligned} f'_{i_1...i_n} &= f\circ(\overline{e}'_{i_1},...,\overline{e}'_{i_n}) \\ &= f\circ(\overline{e}_{j_1}h^{j_1}_{i_1},...,\overline{e}'_{j_n}h^{j_n}_{i_n}) \\ &= h^{j_1}_{i_1}...h^{j_n}_{i_n} f\circ(\overline{e}_{j_1},...,\overline{e}_{j_n}) \\ &= h^{j_1}_{i_1}...h^{j_n}_{i_n} f_{j_1...j_n} \end{aligned}$$

From equation (6.9) the tensor law of transformation of coordinates of polylinear map follows. From equation (6.9) and theorem [4]-2.1.4 it follows that value of the mapping $f\circ(\overline{a}_1,...,\overline{a}_n)$ does not depend from choice of basis. □

Polylinear mapping (5.4) is **symmetric**, if

$$f\circ(a_1,...,a_n) = f\circ(\sigma(a_1),...,\sigma(a_n))$$

for any transposition $\sigma$ of set $\{a_1,...,a_n\}$.

**Theorem 6.3.** *If polyadditive map $f$ is symmetric, then*

(6.10)
$$f_{i_1,...,i_n} = f_{\sigma(i_1),...,\sigma(i_n)}$$

*Proof.* Equation (6.10) follows from equation

$$\begin{aligned} a_1^{i_1}...a_n^{i_n} f_{i_1...i_n} =&f\circ(a_1,...,a_n) \\ =&f\circ(\sigma(a_1),...,\sigma(a_n)) \\ =&a_1^{i_1}...a_n^{i_n} f_{\sigma(i_1)...\sigma(i_n)} \end{aligned}$$

□

Polylinear mapping (5.4) is **skew symmetric**, if

$$f\circ(a_1,...,a_n) = |\sigma| f\circ(\sigma(a_1),...,\sigma(a_n))$$

for any transposition $\sigma$ of set $\{a_1,...,a_n\}$. Here

$$|\sigma| = \begin{cases} 1 & \text{transposition } \sigma \text{ even} \\ -1 & \text{transposition } \sigma \text{ odd} \end{cases}$$

**Theorem 6.4.** *If polylinear map $f$ is skew symmetric, then*

(6.11)
$$f_{i_1,...,i_n} = |\sigma| f_{\sigma(i_1),...,\sigma(i_n)}$$

*Proof.* Equation (6.11) follows from equation

$$\begin{aligned} a_1^{i_1}...a_n^{i_n} f_{i_1...i_n} =&f\circ(a_1,...,a_n) \\ =&|\sigma| f\circ(\sigma(a_1),...,\sigma(a_n)) \\ =&a_1^{i_1}...a_n^{i_n} |\sigma| f_{\sigma(i_1)...\sigma(i_n)} \end{aligned}$$

□

**Theorem 6.5.** *Let $A$ be free finite dimensional associative algebra over the ring $D$. Let $\overline{\overline{I}}$ be basis of algebra $\mathcal{L}(A;A)$. Let $\overline{\overline{e}}$ be the basis of the algebra $A$ over the ring $D$. Let polylinear over ring $D$ mapping $f$ be generated by set of mappings* $(I_{k_1},...,I_{k_n})$. *Coordinates of the mapping* $f$ *and its components relative basis* $\overline{\overline{e}}$ *satisfy to the equation*

$$f_{l_1...l_n}=f^{i_0...i_n}_{t\cdot k_1...k_n}I_{k_1\cdot}{}^{j_1}_{l_1}...I_{k_n\cdot}{}^{j_n}_{l_n}C^{k_1}_{i_0\sigma_t(j_1)}C^{l_1}_{k_1i_1}...B^{k_n}_{l_{n-1}\sigma_t(j_n)}C^{l_n}_{k_ni_n}\overline{e}_{l_n}\tag{6.12}$$

$$f^p_{l_1...l_n}=f^{i_0...i_n}_{t\cdot k_1...k_n}I_{k_1\cdot}{}^{j_1}_{l_1}...I_{k_n\cdot}{}^{j_n}_{l_n}C^{k_1}_{i_0\sigma_t(j_1)}C^{l_1}_{k_1i_1}...C^{k_n}_{l_{n-1}\sigma_t(j_n)}C^{p}_{k_ni_n}\tag{6.13}$$

*Proof.* In equation (6.1), we assume

$$I_{k_i}\circ a_i=\overline{e}_{j_i}I_{k_i\cdot}{}^{j_i}_{l_i}a_i^{l_i}$$

Then equation (6.1) gets form

$$\begin{aligned}f\circ(a_1,...,a_n)&=f^{i_0...i_n}_{t\cdot k_1...k_n}\ \overline{e}_{i_0}\sigma_t(a_1^{l_1}I_{k_1\cdot}{}^{j_1}_{l_1}\overline{e}_{j_1})\overline{e}_{i_1}...\sigma_t(a_n^{l_n}I_{k_n\cdot}{}^{j_n}_{l_n}\overline{e}_{j_n})\overline{e}_{i_n}\\&=a_1^{l_1}...a_n^{l_n}f^{i_0...i_n}_{t\cdot k_1...k_n}I_{k_1\cdot}{}^{j_1}_{l_1}...I_{k_n\cdot}{}^{j_n}_{l_n}\overline{e}_{i_0}\sigma_t(\overline{e}_{j_1})\overline{e}_{i_1}...\sigma_t(\overline{e}_{j_n})\overline{e}_{i_n}\\&=a_1^{l_1}...a_n^{l_n}f^{i_0...i_n}_{t\cdot k_1...k_n}I_{k_1\cdot}{}^{j_1}_{l_1}...I_{k_n\cdot}{}^{j_n}_{l_n}C^{k_1}_{i_0\sigma_t(j_1)}C^{l_1}_{k_1i_1}\\&\quad...C^{k_n}_{l_{n-1}\sigma_t(j_n)}C^{l_n}_{k_ni_n}\overline{e}_{l_n}\end{aligned}\tag{6.14}$$

From equation (6.6) it follows that

$$f\circ(a_1,...,a_n)=\overline{e}_pf^p_{i_1...i_n}a_1^{i_1}...a_n^{i_n}\tag{6.15}$$

Equation (6.12) follows from comparison of equations (6.14) and (6.6). Equation (6.13) follows from comparison of equations (6.14) and (6.15). □

## 8. Index

## 9. Special Symbols and Notations

# Полилинейное отображение свободной алгебры

Александр Клейн

Аннотация. В статье я рассматриваю структуру полилинейного отображения свободной алгебры над коммутативным кольцом.

## Содержание

## 1. Алгебра над кольцом

**Определение 1.1.** Пусть $D$ - коммутативное кольцо. Пусть $A$ - модуль над кольцом $D$.[1] Для заданного билинейного отображения

$$f : A \times A \to A$$

мы определим произведение в $A$

$$ab = f \circ (a, b) \tag{1.1}$$

$A$ - **алгебра над кольцом** $D$, если $A$ - $D$-модуль и в $A$ определена операция произведения (1.1). Если $A$ является свободным $D$-модулем, то $A$ называется **свободной алгеброй над кольцом** $D$. □

Aleks_Kleyn@MailAPS.org.
http://sites.google.com/site/AleksKleyn/.
http://arxiv.org/a/kleyn_a_1.
http://AleksKleyn.blogspot.com/.

[1]Этот подраздел написан на основе раздела [4]-2.2.

Пусть $\overline{\overline{e}}$ - базис свободной алгебры $A$ над кольцом $D$. Если алгебра $A$ имеет единицу, положим $\overline{e}_0$ - единица алгебры $A$.

**Теорема 1.2.** *Пусть $\overline{\overline{e}}$ - базис свободной алгебры $A$ над кольцом $D$. Пусть*

$$a = a^i e_i \quad b = b^i e_i \quad a, b \in A$$

*Произведение $a$, $b$ можно получить согласно правилу*

$$(ab)^k = C_{ij}^k a^i b^j \tag{1.2}$$

*где $C_{ij}^k$ -* **структурные константы алгебры $A$ над кольцом $D$**. *Произведение базисных векторов в алгебре $A$ определено согласно правилу*

$$\overline{e}_i \overline{e}_j = C_{ij}^k \overline{e}_k \tag{1.3}$$

*Доказательство.* Равенство (1.3) является следствием утверждения, что $\overline{\overline{e}}$ является базисом алгебры $A$. Так как произведение в алгебре является билинейным отображением, то произведение $a$ и $b$ можно записать в виде

$$ab = a^i b^j e_i e_j \tag{1.4}$$

Из равенств (1.3), (1.4), следует

$$ab = a^i b^j C_{ij}^k \overline{e}_k \tag{1.5}$$

Так как $\overline{\overline{e}}$ является базисом алгебры $A$, то равенство (1.2) следует из равенства (1.5). □

**Определение 1.3.** Пусть $A_1$ и $A_2$ - алгебры над кольцом $D$. Линейное отображение[2]

$$f : A_1 \to A_2$$

$D$-модуля $A_1$ в $D$-модуль $A_2$ называется **линейным отображением $D$-алгебры $A_1$ в $D$-алгебру $A_2$**. Обозначим $\mathcal{L}(A_1; A_2)$ множество линейных отображений алгебры $A_1$ в алгебру $A_2$. □

**Определение 1.4.** Пусть $D$ - коммутативное ассоциативное кольцо. Пусть $A_1$, ..., $A_n$ - $D$-алгебры и $S$ - $D$-модуль. Мы будем называть отображение

$$f : A_1 \times ... \times A_n \to S$$

**полилинейным отображением алгебр** $A_1$, ..., $A_n$ в модуль $S$, если

$$f \circ (a_1, ..., a_i + b_i, ..., a_n) = f \circ (a_1, ..., a_i, ..., a_n) + f \circ (a_1, ..., b_i, ..., a_n)$$
$$f \circ (a_1, ..., pa_i, ..., a_n) = pf \circ (a_1, ..., a_i, ..., a_n)$$

$$1 \le i \le n \quad a_i, b_i \in A_i \quad p \in D$$

Обозначим $\mathcal{L}(A_1, ..., A_n; S)$ множество полилинейных отображений алгебр $A_1$, ..., $A_n$ в модуль $S$. Обозначим $\mathcal{L}(A^n; S)$ множество $n$-линейных отображений алгебры $A$ ($A_1 = ... = A_n = A$) в модуль $S$. □

[2]Этот подраздел написан на основе раздела [4]-2.3.

**Теорема 1.5.** *Пусть $D$ - коммутативное ассоциативное кольцо. Пусть $A_1$, ..., $A_n$ - $D$-алгебры и $S$ - $D$-модуль. Пусть отображения*

$$f : A_1 \times ... \times A_n \to S$$
$$g : A_1 \times ... \times A_n \to S$$

*являются полилинейными отображениями. Тогда отображение $f+g$, определённое равенством*

$$(f+g) \circ (a_1, ..., a_n) = f \circ (a_1, ..., a_n) + g \circ (a_1, ..., a_n)$$

*также является полилинейным.*

*Доказательство.* Утверждение теоремы следует из цепочек равенств

$$\begin{aligned}
&(f+g) \circ (x_1, ..., x_i + y_i, ..., x_n)\\
=&f \circ (x_1, ..., x_i + y_i, ..., x_n) + g \circ (x_1, ..., x_i + y_i, ..., x_n)\\
=&f \circ (x_1, ..., x_i, ..., x_n) + f \circ (x_1, ..., y_i, ..., x_n)\\
+&g \circ (x_1, ..., x_i, ..., x_n) + g \circ (x_1, ..., y_i, ..., x_n)\\
=&(f+g) \circ (x_1, ..., x_i, ..., x_n) + (f+g) \circ (x_1, ..., y_i, ..., x_n)
\end{aligned}$$

$$\begin{aligned}
&(f+g) \circ (x_1, ..., px_i, ..., x_n)\\
=&f \circ (x_1, ..., px_i, ..., x_n) + g \circ (x_1, ..., px_i, ..., x_n)\\
=&pf \circ (x_1, ..., x_i, ..., x_n) + pg \circ (x_1, ..., x_i, ..., x_n)\\
=&p(f \circ (x_1, ..., x_i, ..., x_n) + g \circ (x_1, ..., x_i, ..., x_n))\\
=&p(f+g) \circ (x_1, ..., x_i, ..., x_n)
\end{aligned}$$

□

**Следствие 1.6.** *Рассмотрим алгебру $A_1$ и алгебру $A_2$. Пусть отображения*

$$f : A_1 \to A_2$$
$$g : A_1 \to A_2$$

*являются линейными отображениями. Тогда отображение $f+g$, определённое равенством*

$$(f+g) \circ a = f \circ a + g \circ a$$

*также является линейным.* □

**Теорема 1.7.** *Пусть $D$ - коммутативное ассоциативное кольцо. Пусть $A_1$, ..., $A_n$ - $D$-алгебры и $S$ - $D$-модуль. Пусть отображение*

$$f : A_1 \times ... \times A_n \to S$$

*является полилинейным отображением. Тогда отображение $pf$, $p \in D$, определённое равенством*

$$(pf) \circ x = p\ f \circ x$$

*также является полилинейными. При этом выполняется равенство*

$$p(qf) = (pq)f$$
$$(p+q)f = pf + qf$$

*Доказательство.* Утверждение теоремы следует из цепочек равенств

$$\begin{aligned}(pf)\circ(x_1,...,x_i+y_i,...,x_n)=&p\ f\circ(x_1,...,x_i+y_i,...,x_n)\\=&p\ (f\circ(x_1,...,x_i,...,x_n)+f\circ(x_1,...,y_i,...,x_n))\\=&p\ f\circ(x_1,...,x_i,...,x_n)+p\ f\circ(x_1,...,y_i,...,x_n)\\=&(pf)\circ(x_1,...,x_i,...,x_n)+(pf)\circ(x_1,...,y_i,...,x_n)\\(pf)\circ(x_1,...,qx_i,...,x_n)=&p\ f\circ(x_1,...,qx_i,...,x_n)=pq\ f\circ(x_1,...,x_i,...,x_n)\\=&qp\ f\circ(x_1,...,x_n)=q\ (pf)\circ(x_1,...,x_n)\\(p(qf))\circ(x_1,...,x_n)=&p\ (qf)\circ(x_1,...,x_n)=p\ (q\ f\circ(x_1,...,x_n))\\=&(pq)\ f\circ(x_1,...,x_n)=((pq)f)\circ(x_1,...,x_n)\\((p+q)f)\circ(x_1,...,x_n)=&(p+q)\ f\circ(x_1,...,x_n)\\=&p\ f\circ(x_1,...,x_n)+q\ f\circ(x_1,...,x_n)\\=&(pf)\circ(x_1,...,x_n)+(qf)\circ(x_1,...,x_n)\end{aligned}$$

□

**Следствие 1.8.** *Рассмотрим алгебру $A_1$ и алгебру $A_2$. Пусть отображение*

$$f:A_1\to A_2$$

*является линейным отображением. Тогда отображение $pf$, $p\in D$, определённое равенством*

$$(pf)\circ x=p\ f\circ x$$

*также является линейным. При этом выполняется равенство*

$$p(qf)=(pq)f$$

$$(p+q)f=pf+qf$$

□

**Теорема 1.9.** *Пусть $D$ - коммутативное ассоциативное кольцо. Пусть $A_1$, ..., $A_n$ - $D$-алгебры и $S$ - $D$-модуль. Множество $\mathcal{L}(A_1,...,A_n;S)$ является $D$-модулем.*

*Доказательство.* Теорема 1.5 определяет сумму полилинейных отображений в $D$-модуль $S$. Пусть $f$, $g$, $h\in\mathcal{L}(A_1,...,A_2;S)$. Для любого $a=(a_1,...,a_n)$, $a_1\in A_1$, ..., $a_n\in A_n$,

$$\begin{aligned}(f+g)\circ a=&f\circ a+g\circ a=g\circ a+f\circ a\\=&(g+f)\circ a\\((f+g)+h)\circ a=&(f+g)\circ a+h\circ a=(f\circ a+g\circ a)+h\circ a\\=&f\circ a+(g\circ a+h\circ a)=f\circ a+(g+h)\circ a\\=&(f+(g+h))\circ a\end{aligned}$$

Следовательно, сумма полилинейных отображений коммутативна и ассоциативна.

Отображение $z$, определённое равенством

$$z\circ a=0$$

является нулём операции сложения, так как

$$(z+f)\circ a=z\circ a+f\circ a=0+f\circ a=f\circ a$$

Для заданного отображения $f$ отображение $g$, определённое равенством

$$g \circ a = -f \circ a$$

удовлетворяет равенству

$$f + g = z$$

так как

$$(f+g) \circ a = f \circ a + g \circ a = f \circ a - f \circ a = 0$$

Следовательно, множество $\mathcal{L}(A_1; A_2)$ является абелевой группой.

Из теоремы 1.7 следует, что определенно представление кольца $D$ в абелевой группе $\mathcal{L}(A_1, ..., A_n; S)$. Так как кольцо $D$ имеет единицу, то согласно теореме [4]-2.1.1 указанное представление эффективно. □

**Следствие 1.10.** *Пусть $D$ - коммутативное кольцо с единицей. Рассмотрим $D$-алгебру $A_1$ и $D$-алгебру $A_2$. Множество $\mathcal{L}(A_1; A_2)$ является $D$-модулем.* □

**Теорема 1.11.** *Пусть $A$ - алгебра над коммутативным кольцом $D$. Модуль $\mathcal{L}(A; A)$, оснащённый произведением*

$$\circ : (g, f) \in \mathcal{L}(A; A) \times \mathcal{L}(A; A) \to g \circ f \in \mathcal{L}(A; A) \tag{1.6}$$

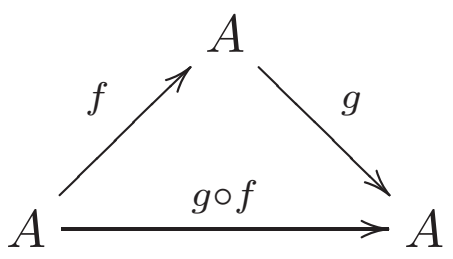

*является алгеброй над $D$.*

*Доказательство.* Смотри доказательство теоремы [4]-2.4.5. □

## 2. Тензорное произведение алгебр

**Определение 2.1.** Пусть $A_1$, ..., $A_n$ - свободные алгебры над коммутативным кольцом $D$.[3] Рассмотрим категорию $\mathcal{A}$ объектами которой являются полилинейные над коммутативным кольцом $D$ отображения

$$f : A_1 \times ... \times A_n \longrightarrow S_1 \qquad g : A_1 \times ... \times A_n \longrightarrow S_2$$

где $S_1$, $S_2$ - модули над кольцом $D$. Мы определим морфизм $f \to g$ как линейное над коммутативным кольцом $D$ отображение $h : S_1 \to S_2$, для которого коммутативна диаграмма

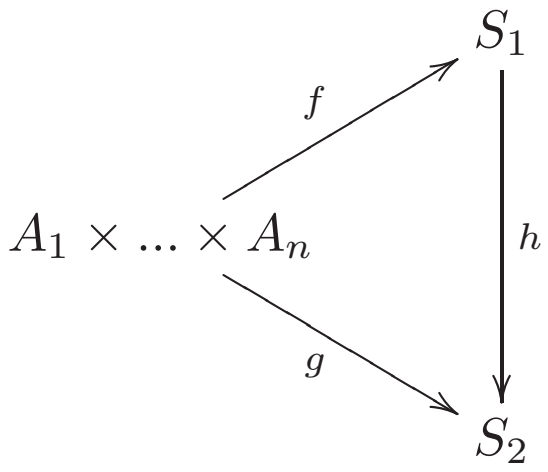

Универсальный объект $A_1 \otimes ... \otimes A_n$ категории $\mathcal{A}$ называется **тензорным произведением алгебр** $A_1$, ..., $A_n$. □

[3] Я определяю тензорное произведение алгебр по аналогии с определением в [1], с. 456 - 458. Этот подраздел написан на основе раздела [4]-2.5.

**Определение 2.2.** Тензорное произведение

$$A^{\otimes n} = A_1 \otimes ... \otimes A_n \quad A_1 = ... = A_n = A$$

называется **тензорной степенью алгебры** $A$. □

**Теорема 2.3.** *Тензорное произведение алгебр существует.*

*Доказательство.* Пусть $M$ - модуль над кольцом $D$, порождённый произведением $A_1 \times ... \times A_n$ алгебр $A_1$, ..., $A_n$. Инъекция

$$i : A_1 \times ... \times A_n \longrightarrow M$$

определена по правилу

(2.1) $$i \circ (d_1, ..., d_n) = (d_1, ..., d_n)$$

Пусть $N \subset M$ - подмодуль, порождённый элементами вида

(2.2) $$(d_1, ..., d_i + c_i, ..., d_n) - (d_1, ..., d_i, ..., d_n) - (d_1, ..., c_i, ..., d_n)$$

(2.3) $$(d_1, ..., ad_i, ..., d_n) - a(d_1, ..., d_i, ..., d_n)$$

где $d_i \in A_i$, $c_i \in A_i$, $a \in D$. Пусть

$$j : M \to M/N$$

каноническое отображение на фактормодуль. Рассмотрим коммутативную диаграмму

(2.4)

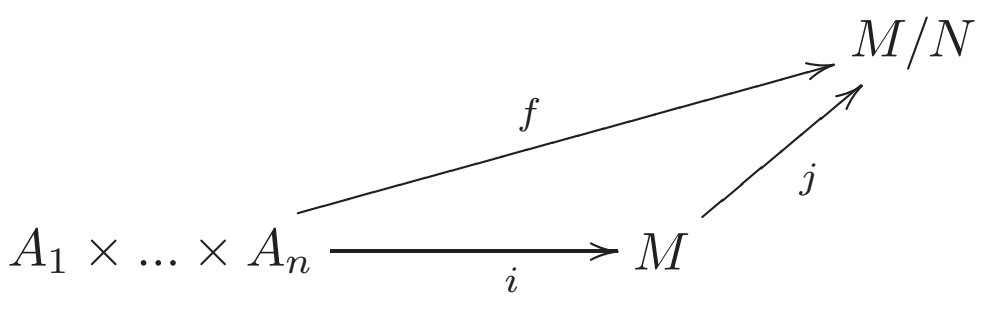

Поскольку элементы (2.2) и (2.3) принадлежат ядру линейного отображения $j$, то из равенства (2.1) следует

(2.5) $$f \circ (d_1, ..., d_i + c_i, ..., d_n) = f \circ (d_1, ..., d_i, ..., d_n) + f \circ (d_1, ..., c_i, ..., d_n)$$

(2.6) $$f \circ (d_1, ..., ad_i, ..., d_n) = a\ f \circ (d_1, ..., d_i, ..., d_n)$$

Из равенств (2.5) и (2.6) следует, что отображение $f$ полилинейно над кольцом $D$. Поскольку $M$ - модуль с базисом $A_1 \times ... \times A_n$, то, согласно теореме [1]-1 на с. 104, для любого модуля $V$ и любого полилинейного над $D$ отображения

$$g : A_1 \times ... \times A_n \longrightarrow V$$

существует единственный гомоморфизм $k : M \to V$, для которого коммутативна следующая диаграмма

(2.7)

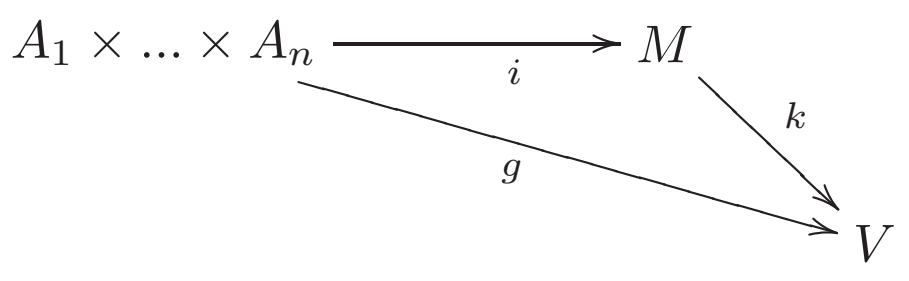

Так как $g$ - полилинейно над $D$, то $\ker k \subseteq N$. Согласно утверждению на с. [1]-94, отображение $j$ универсально в категории гомоморфизмов векторного пространства $M$, ядро которых содержит $N$. Следовательно, определён гомоморфизм

$$h : M/N \to V$$

для которого коммутативна диаграмма

(2.8)

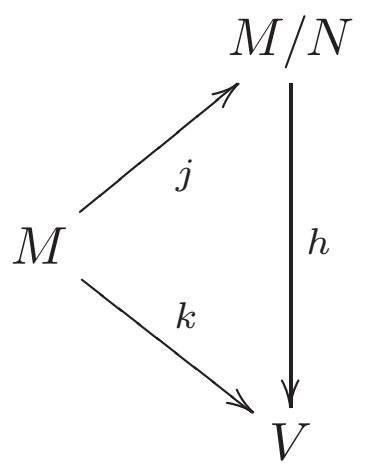

Объединяя диаграммы (2.4), (2.7), (2.8), получим коммутативную диаграмму

(2.9)

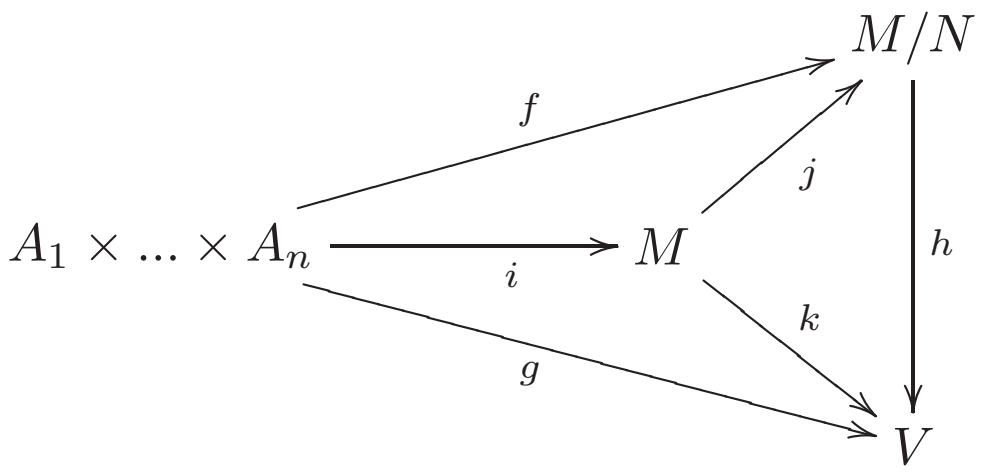

Так как Im$f$ порождает $M/N$, то отображение $h$ однозначно определено. □

Согласно доказательству теоремы 2.3

$$A_1 \otimes ... \otimes A_n = M/N$$

Для $d_i \in A_i$ будем записывать

$$j \circ (d_1, ..., d_n) = d_1 \otimes ... \otimes d_n \tag{2.10}$$

**Теорема 2.4.** *Пусть $A_1$, ..., $A_n$ - алгебры над коммутативным кольцом $D$. Пусть*

$$f : A_1 \times ... \times A_n \to A_1 \otimes ... \otimes A_n$$

*полилинейное отображение, определённое равенством*

$$f \circ (d_1, ..., d_n) = d_1 \otimes ... \otimes d_n \tag{2.11}$$

*Пусть*

$$g : A_1 \times ... \times A_n \to V$$

*полилинейное отображение в $D$-модуль $V$. Существует $D$-линейное отображение*

$$h : A_1 \otimes ... \otimes A_n \to V$$

*такое, что диаграмма*

(2.12)

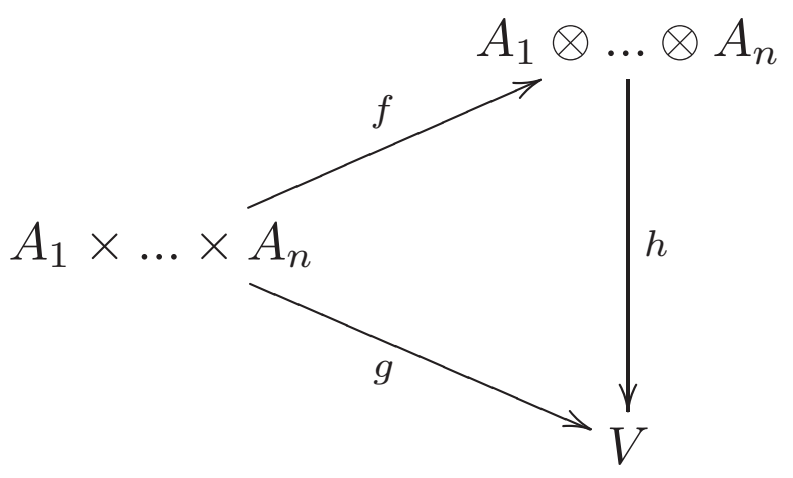

*коммутативна.*

*Доказательство.* Равенство (2.11) следует из равенств (2.1) и (2.10). Существование отображения $h$ следует из определения 2.1 и построений, выполненных при доказательстве теоремы 2.3. □

Равенства (2.5) и (2.6) можно записать в виде

$$\begin{aligned}&a_1\otimes ...\otimes(a_i+b_i)\otimes ...\otimes a_n\\=&a_1\otimes ...\otimes a_i\otimes ...\otimes a_n+a_1\otimes ...\otimes b_i\otimes ...\otimes a_n\end{aligned}\tag{2.13}$$

$$a_1\otimes ...\otimes(ca_i)\otimes ...\otimes a_n=c(a_1\otimes ...\otimes a_i\otimes ...\otimes a_n)\tag{2.14}$$

$$a_i\in A_i\quad b_i\in A_i\quad c\in D$$

**Теорема 2.5.** *Пусть $A$ - алгебра над коммутативным кольцом $D$. Существует линейное отображение*

$$h:a\otimes b\in A\otimes A\to ab\in A$$

*Доказательство.* Теорема является следствием теоремы 2.4 и определения 1.1. □

**Теорема 2.6.** *Тензорное произведение $A_1\otimes ...\otimes A_n$ свободных конечномерных алгебр $A_1$, ..., $A_n$ над коммутативным кольцом $D$ является свободной конечномерной алгеброй.*

*Пусть $\overline{\overline{e}}_i$ - базис алгебры $A_i$ над кольцом $D$. Произвольный тензор $a\in A_1\otimes ...\otimes A_n$ можно представить в виде*

$$a=a^{i_1...i_n}\overline{e}_{1\cdot i_1}\otimes ...\otimes\overline{e}_{n\cdot i_n}\tag{2.15}$$

*Мы будем называть выражение $a^{i_1...i_n}$* **стандартной компонентой тензора***.*

*Доказательство.* Смотри доказательство теоремы [4]-2.5.6. □

## 3. Линейное отображение в ассоциативную алгебру

**Теорема 3.1.** *Рассмотрим $D$-алгебры $A_1$ и $A_2$. Для заданного отображения $f\in\mathcal L(A_1;A_2)$ существует линейное отображение*

$$h:A_2\otimes A_2\to\mathcal L(A_1;A_2)$$

*определённое равенством*

$$(a\otimes b)\circ f=afb\tag{3.1}$$

*Доказательство.* Смотри доказательство теорем [4]-2.6.1 и [4]-2.6.2. □

**Теорема 3.2.** *Рассмотрим $D$-алгебры $A_1$ и $A_2$. Определим произведение в алгебре $A_2\otimes A_2$ согласно правилу*

$$(c\otimes d)\circ(a\otimes b)=(ca)\otimes(bd)\tag{3.2}$$

*Линейное отображение*

$$h:A_2\otimes A_2\to{}^*\mathcal L(A_1;A_2)\tag{3.3}$$

*определённое равенством*

$$(a\otimes b)\circ f=afb\quad a,b\in A_2\quad f\in\mathcal L(A_1;A_2)\tag{3.4}$$

*является представлением*[4] *алгебры* $A_2 \otimes A_2$ *в модуле* $\mathcal{L}(A_1; A_2)$.

*Доказательство.* Смотри доказательство теоремы [4]-2.6.3 □

**Теорема 3.3.** *Рассмотрим* $D$*-алгебру* $A$. *Определим произведение в алгебре* $A \otimes A$ *согласно правилу* (3.2). *Представление*

$$h : A \otimes A \to {}^*\mathcal{L}(A; A) \tag{3.5}$$

*алгебры* $A \otimes A$ *в модуле* $\mathcal{L}(A; A)$, *определённое равенством*

$$(a \otimes b) \circ f = afb \quad a, b \in A \quad f \in \mathcal{L}(A; A) \tag{3.6}$$

*позволяет отождествить тензор* $d \in A \otimes A$ *с отображением* $d \circ \delta \in \mathcal{L}(A; A)$, *где* $\delta \in \mathcal{L}(A; A)$ *- тождественное отображение.*

*Доказательство.* Смотри доказательство теоремы [4]-2.6.4 □

Из теоремы 3.3 следует, что отображение (3.4) можно рассматривать как произведение отображений $a \otimes b$ и $f$.

**Теорема 3.4.** *Рассмотрим* $D$*-алгебру* $A_1$ *и ассоциативную* $D$*-алгебру* $A_2$. *Рассмотрим представление алгебры* $A_2 \otimes A_2$ *в модуле* $\mathcal{L}(A_1; A_2)$. *Отображение*

$$h : A_1 \to A_2$$

*порождённое отображением*

$$f : A_1 \to A_2$$

*имеет вид*

$$h = (a_{s\cdot 0} \otimes a_{s\cdot 1}) \circ f = a_{s\cdot 0} f a_{s\cdot 1} \tag{3.7}$$

*Доказательство.* Смотри доказательство теоремы [4]-2.6.6 □

**Теорема 3.5.** *Пусть* $A$ *- свободная конечно мерная ассоциативная алгебра над полем* $D$. *Представление алгебры* $A \otimes A$ *в алгебре* $\mathcal{L}(A; A)$ *имеет конечный* **базис** $\overline{\overline{I}}$.

(1) *Линейное отображение* $f \in \mathcal{L}(A; A)$ *имеет вид*

$$f = (a_{k\cdot s_k\cdot 0} \otimes a_{k\cdot s_k\cdot 1}) \circ I_k = \sum_{k} a_{k\cdot s_k\cdot 0} I_k a_{k\cdot s_k\cdot 1} \tag{3.8}$$

(2) *Его стандартное представление имеет вид*

$$f = a^{k\cdot ij}(\overline{e}_i \otimes \overline{e}_j) \circ I_k = a^{k\cdot ij}\overline{e}_i I_k \overline{e}_j \tag{3.9}$$

*Доказательство.* Смотри доказательство теоремы [4]-2.7.5 □

[4] Определение представления $\Omega$-алгебры дано в определении [3]-2.1.4.

## 4. Линейное отображение в неассоциативную алгебру

Так как произведение неассоциативно, мы можем предположить, что действие $a, b \in A$ на отображение $f$ может быть представленно либо в виде $a(fb)$, либо в виде $(af)b$.

**Теорема 4.1.** *Пусть $\overline{\overline{e}}_1$ - базис свободной конечно мерной $D$-алгебры $A_1$. Пусть $\overline{\overline{e}}_2$ - базис свободной конечно мерной неассоциативной $D$-алгебры $A_2$. Пусть $C_2{}_{\cdot}{}^{p}_{kl}$ - структурные константы алгебры $A_2$. Пусть отображение*

$$g = a \circ f \tag{4.1}$$

*порождённое отображением $f \in (A_1; A_2)$ посредством тензора $a \in A_2 \otimes A_2$, имеет стандартное представление*

$$g = a^{ij}(\overline{e}_i \otimes \overline{e}_j) \circ f = a^{ij}(\overline{e}_i f)\overline{e}_j \tag{4.2}$$

*Координаты отображения (4.1) и его стандартные компоненты связаны равенством*

$$g^k_l = f^m_l g^{ij} C_2{}_{\cdot}{}^{p}_{im} C_2{}_{\cdot}{}^{k}_{pj} \tag{4.3}$$

*Доказательство.* Смотри доказательство теоремы [4]-2.8.3 □

**Теорема 4.2.** *Пусть $A$ - свободная конечно мерная неассоциативная алгебра над кольцом $D$. Представление алгебры $A \otimes A$ в алгебре $\mathcal{L}(A; A)$ имеет конечный базис $\overline{\overline{I}}$.*

(1) *Линейное отображение $f \in \mathcal{L}(A; A)$ имеет вид*

$$f = (a_{k \cdot s_k \cdot 0} \otimes a_{k \cdot s_k \cdot 1}) \circ I_k = (a_{k \cdot s_k \cdot 0} I_k) a_{k \cdot s_k \cdot 1} \tag{4.4}$$

(2) *Его стандартное представление имеет вид*

$$f = a^{k \cdot ij}(\overline{e}_i \otimes \overline{e}_j) \circ I_k = a^{k \cdot ij}(\overline{e}_i I_k)\overline{e}_j \tag{4.5}$$

*Доказательство.* Смотри доказательство теоремы [4]-2.8.4 □

## 5. Полилинейное отображение в ассоциативную алгебру

**Теорема 5.1.** *Пусть $A_1$, ..., $A_n$, $A$ - ассоциативные $D$-алгебры. Пусть*

$$f_i \in \mathcal{L}(A_i; A) \quad i = 1, ..., n$$

$$a_j \in A \quad j = 0, ..., n$$

*Для заданной перестановки $\sigma$ $n$ переменных отображение*

$$\begin{aligned} &((a_0, ..., a_n) \circ \sigma(f_1, ..., f_n)) \circ (x_1, ..., x_n) \\ =& (a_0 \sigma(f_1) a_1 ... a_{n-1} \sigma(f_n) a_n) \circ (x_1, ..., x_n) \\ =& a_0 \sigma(f_1 \circ x_1) a_1 ... a_{n-1} \sigma(f_n \circ x_n) a_n \end{aligned} \tag{5.1}$$

*является $n$-линейным отображением в алгебру $A$.*

*Доказательство.* Утверждение теоремы следует из цепочек равенств

$$\begin{aligned}
&((a_0,...,a_n,\sigma)\circ(f_1,...,f_n))\circ(x_1,...,x_i+y_i,...,x_n)\\
=&a_0\sigma(f_1\circ x_1)a_1...\sigma(f_i\circ(x_i+y_i))...a_{n-1}\sigma(f_n\circ x_n)a_n\\
=&a_0\sigma(f_1\circ x_1)a_1...\sigma(f_i\circ x_i+f_i\circ y_i)...a_{n-1}\sigma(f_n\circ x_n)a_n\\
=&a_0\sigma(f_1\circ x_1)a_1...\sigma(f_i\circ x_i)...a_{n-1}\sigma(f_n\circ x_n)a_n\\
+&a_0\sigma(f_1\circ x_1)a_1...\sigma(f_i\circ y_i)...a_{n-1}\sigma(f_n\circ x_n)a_n\\
=&((a_0,...,a_n,\sigma)\circ(f_1,...,f_n))\circ(x_1,...,x_i,...,x_n)\\
+&((a_0,...,a_n,\sigma)\circ(f_1,...,f_n))\circ(x_1,...,y_i,...,x_n)
\end{aligned}$$

$$\begin{aligned}
&((a_0,...,a_n,\sigma)\circ(f_1,...,f_n))\circ(x_1,...,px_i,...,x_n)\\
=&a_0\sigma(f_1\circ x_1)a_1...\sigma(f_i\circ(px_i))...a_{n-1}\sigma(f_n\circ x_n)a_n\\
=&a_0\sigma(f_1\circ x_1)a_1...\sigma(p(f_i\circ x_i))...a_{n-1}\sigma(f_n\circ x_n)a_n\\
=&p(a_0\sigma(f_1\circ x_1)a_1...\sigma(f_i\circ x_i)...a_{n-1}\sigma(f_n\circ x_n)a_n)\\
=&p(((a_0,...,a_n,\sigma)\circ(f_1,...,f_n))\circ(x_1,...,x_i,...,x_n))
\end{aligned}$$

□

В равенстве (5.1), также как и в других выражениях полилинейного отображения, принято соглашение, что отображение $f_i$ имеет своим аргументом переменную $x_i$.

**Теорема 5.2.** *Пусть $A_1$, ..., $A_n$, $A$ - ассоциативные $D$-алгебры. Для заданного семейства отображений*

$$f_i\in\mathcal{L}(A_i;A)\quad i=1,...,n$$

*отображение*

$$h:A^{n+1}\to\mathcal{L}(A_1,...,A_n;A)$$

*определённое равенством*

$$(a_0,...,a_n,\sigma)\circ(f_1,...,f_n)=a_0\sigma(f_1)a_1...a_{n-1}\sigma(f_n)a_n$$

*является $n+1$-линейным отображением в $D$-модуль $\mathcal{L}(A_1,...,A_n;A)$.*

*Доказательство.* Утверждение теоремы следует из цепочек равенств

$$\begin{aligned}
&((a_0,...,a_i+b_i,...a_n,\sigma)\circ(f_1,...,f_n))\circ(x_1,...,x_n)\\
=&a_0\sigma(f_1\circ x_1)a_1...(a_i+b_i)...a_{n-1}\sigma(f_n\circ x_n)a_n\\
=&a_0\sigma(f_1\circ x_1)a_1...a_i...a_{n-1}\sigma(f_n\circ x_n)a_n+a_0\sigma(f_1\circ x_1)a_1...b_i...a_{n-1}\sigma(f_n\circ x_n)a_n\\
=&((a_0,...,a_i,...,a_n,\sigma)\circ(f_1,...,f_n))\circ(x_1,...,x_n)\\
+&((a_0,...,b_i,...,a_n,\sigma)\circ(f_1,...,f_n))\circ(x_1,...,x_n)\\
=&((a_0,...,a_i,...,a_n,\sigma)\circ(f_1,...,f_n)+(a_0,...,b_i,...,a_n,\sigma)\circ(f_1,...,f_n))\circ(x_1,...,x_n)
\end{aligned}$$

$$\begin{aligned}
&((a_0,...,pa_i,...a_n,\sigma)\circ(f_1,...,f_n))\circ(x_1,...,x_n)\\
=&a_0\sigma(f_1\circ x_1)a_1...pa_i...a_{n-1}\sigma(f_n\circ x_n)a_n\\
=&p(a_0\sigma(f_1\circ x_1)a_1...a_i...a_{n-1}\sigma(f_n\circ x_n)a_n)\\
=&p(((a_0,...,a_i,...,a_n,\sigma)\circ(f_1,...,f_n))\circ(x_1,...,x_n))\\
=&(p((a_0,...,a_i,...,a_n,\sigma)\circ(f_1,...,f_n)))\circ(x_1,...,x_n)
\end{aligned}$$

□

**Теорема 5.3.** *Пусть $A_1$, ..., $A_n$, $A$ - ассоциативные $D$-алгебры. Для заданного семейства отображений*

$$f_i \in \mathcal{L}(A_i; A) \quad i = 1, ..., n$$

*существует линейное отображение*

$$h : A^{\otimes n+1} \times S_n \to \mathcal{L}(A_1, ..., A_n; A)$$

*определённое равенством*

$$\begin{aligned}(a_0 \otimes ... \otimes a_n, \sigma) \circ (f_1, ..., f_n) &= (a_0, ..., a_n, \sigma) \circ (f_1, ..., f_n)\\ &= a_0\sigma(f_1)a_1...a_{n-1}\sigma(f_n)a_n\end{aligned} \tag{5.2}$$

*Доказательство.* Утверждение теоремы является следствием теорем 2.4, 5.2. □

**Теорема 5.4.** *Пусть $A_1$, ..., $A_n$, $A$ - ассоциативные $D$-алгебры. Для заданного тензора $a \in A^{\otimes n+1}$ и заданной перестановки $\sigma \in S_n$ отображение*

$$h : \prod_{i=1}^{n} \mathcal{L}(A_i; A) \to \mathcal{L}(A_1, ..., A_n; A)$$

*определённое равенством*

$$(a_0 \otimes ... \otimes a_n, \sigma) \circ (f_1, ..., f_n) = a_0\sigma(f_1)a_1...a_{n-1}\sigma(f_n)a_n$$

*является $n$-линейным отображением в $D$-модуль $\mathcal{L}(A_1, ..., A_n; A)$.*

*Доказательство.* Утверждение теоремы следует из цепочек равенств

$$\begin{aligned}
&((a_0\otimes ...\otimes a_n,\sigma)\circ(f_1,...,f_i+g_i,...,f_n))\circ(x_1,...,x_n)\\
=&(a_0\sigma(f_1)a_1...\sigma(f_i+g_i)...a_{n-1}\sigma(f_n)a_n)\circ(x_1,...,x_n)\\
=&a_0\sigma(f_1\circ x_1)a_1...\sigma((f_i+g_i)\circ x_i)...a_{n-1}\sigma(f_n\circ x_n)a_n\\
=&a_0\sigma(f_1\circ x_1)a_1...\sigma(f_i\circ x_i+g_i\circ x_i)...a_{n-1}\sigma(f_n\circ x_n)a_n\\
=&a_0\sigma(f_1\circ x_1)a_1...\sigma(f_i\circ x_i)...a_{n-1}\sigma(f_n\circ x_n)a_n\\
+&a_0\sigma(f_1\circ x_1)a_1...\sigma(g_i\circ x_i)...a_{n-1}\sigma(f_n\circ x_n)a_n\\
=&(a_0\sigma(f_1)a_1...\sigma(f_i)...a_{n-1}\sigma(f_n)a_n)\circ(x_1,...,x_n)\\
+&(a_0\sigma(f_1)a_1...\sigma(g_i)...a_{n-1}\sigma(f_n)a_n)\circ(x_1,...,x_n)\\
=&((a_0\otimes ...\otimes a_n,\sigma)\circ(f_1,...,f_i,...,f_n))\circ(x_1,...,x_n)\\
+&((a_0\otimes ...\otimes a_n,\sigma)\circ(f_1,...,g_i,...,f_n))\circ(x_1,...,x_n)\\
=&((a_0\otimes ...\otimes a_n,\sigma)\circ(f_1,...,f_i,...,f_n)\\
+&(a_0\otimes ...\otimes a_n,\sigma)\circ(f_1,...,g_i,...,f_n))\circ(x_1,...,x_n)
\end{aligned}$$

$$\begin{aligned}
&((a_0\otimes ...\otimes a_n,\sigma)\circ(f_1,...,pf_i,...,f_n))\circ(x_1,...,x_n)\\
=&(a_0\sigma(f_1)a_1,...\sigma(pf_i)...a_{n-1}\sigma(f_n)a_n)\circ(x_1,...,x_n)\\
=&a_0\sigma(f_1\circ x_1)a_1,...\sigma((pf_i)\circ x_i)...a_{n-1}\sigma(f_n\circ x_n)a_n\\
=&a_0\sigma(f_1\circ x_1)a_1,...\sigma(p(f_i\circ x_i))...a_{n-1}\sigma(f_n\circ x_n)a_n\\
=&p(a_0\sigma(f_1\circ x_1)a_1,...\sigma(f_i\circ x_i)...a_{n-1}\sigma(f_n\circ x_n)a_n)\\
=&p(((a_0\otimes ...\otimes a_n,\sigma)\circ(f_1,...,f_i,...,f_n))\circ(x_1,...,x_n))\\
=&(p((a_0\otimes ...\otimes a_n,\sigma)\circ(f_1,...,f_i,...,f_n)))\circ(x_1,...,x_n)
\end{aligned}$$

□

**Теорема 5.5.** *Пусть $A_1$, ..., $A_n$, $A$ - ассоциативные $D$-алгебры. Для заданного тензора $a\in A^{\otimes n+1}$ и заданной перестановки $\sigma\in S_n$ существует линейное отображение*

$$h:\mathcal{L}(A_1;A)\otimes ...\otimes\mathcal{L}(A_n;A)\to\mathcal{L}(A_1,...,A_n;A)$$

*определённое равенством*

$$(a_0\otimes ...\otimes a_n,\sigma)\circ(f_1\otimes ...\otimes f_n)=(a_0\otimes ...\otimes a_n,\sigma)\circ(f_1,...,f_n)\tag{5.3}$$

*Доказательство.* Утверждение теоремы является следствием теорем 2.4, 5.4. □

**Теорема 5.6.** *Пусть $A$ - ассоциативная $D$-алгебра. Полилинейное отображение*

$$f:A^n\to A, a=f\circ(a_1,...,a_n)\tag{5.4}$$

*имеет вид*

$$a=f^n_{s\cdot 0}\ \sigma_s(I_{1\cdot s}\circ a_1)\ f^n_{s\cdot 1}\ ...\ \sigma_s(I_{n\cdot s}\circ a_n)\ f^n_{s\cdot n}\tag{5.5}$$

*где $\sigma_s$ - перестановка множества переменных $\{a_1,...,a_n\}$*

$$\sigma_s=\begin{pmatrix}a_1&...&a_n\\\sigma_s(a_1)&...&\sigma_s(a_n)\end{pmatrix}$$

*Доказательство.* Мы докажем утверждение индукцией по $n$.

При $n = 1$ доказываемое утверждение является следствием утверждения (1) теоремы 3.5. При этом мы можем отождествить[5]

$$f^1_{s\cdot p} = f_{s\cdot p} \quad p = 0, 1$$

Допустим, что утверждение теоремы справедливо при $n = k-1$. Тогда отображение (5.4) можно представить в виде

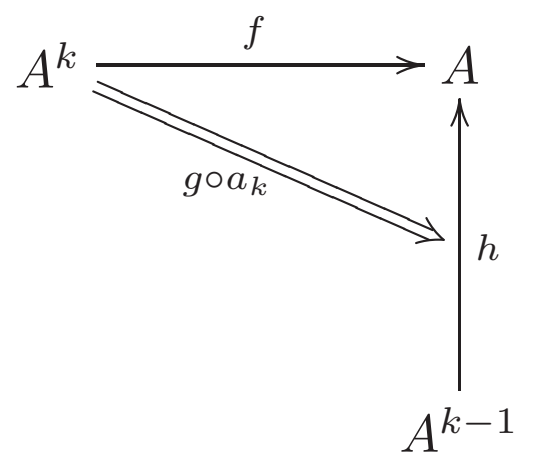

$$a = f \circ (a_1, ..., a_k) = (g \circ a_k) \circ (a_1, ..., a_{k-1})$$

Согласно предположению индукции полилинейное отображение $h$ имеет вид

$$a = h^{k-1}_{t\cdot 0}\ \sigma_t(I_{1\cdot t} \circ a_1)\ h^{k-1}_{t\cdot 1}\ ...\ \sigma_t(I_{k-1\cdot t} \circ a_{k-1})\ h^{k-1}_{t\cdot k-1}$$

Согласно построению $h = g \circ a_k$. Следовательно, выражения $h_{t\cdot p}$ являются функциями $a_k$. Поскольку $g \circ a_k$ - линейная функция $a_k$, то только одно выражение $h_{t\cdot p}$ является линейной функцией переменной $a_k$, и остальные выражения $h_{t\cdot q}$ не зависят от $a_k$.

Не нарушая общности, положим $p = 0$. Согласно равенству (3.7) для заданного $t$

$$h^{k-1}_{t\cdot 0} = g_{tr\cdot 0}\ I_{k\cdot r} \circ a_k\ g_{tr\cdot 1}$$

Положим $s = tr$ и определим перестановку $\sigma_s$ согласно правилу

$$\sigma_s = \sigma_{tr} = \begin{pmatrix} a_k & a_1 & ... & a_{k-1} \\ a_k & \sigma_t(a_1) & ... & \sigma_t(a_{k-1}) \end{pmatrix}$$

Положим

$$f^k_{tr\cdot q+1} = h^{k-1}_{t\cdot q} \quad q = 1, ..., k-1$$
$$f^k_{tr\cdot q} = g_{tr\cdot q} \quad q = 0, 1$$

Мы доказали шаг индукции. □

**Определение 5.7.** Выражение $f^n_{s\cdot p}$ в равенстве (5.5) называется **компонентой полилинейного отображения** $f$. □

**Теорема 5.8.** *Рассмотрим $D$-алгебру $A$. Линейное отображение*

$$h : A^{\otimes n+1} \times S_n \to {}^*\mathcal{L}(A^n; A)$$

[5] В представлении (5.5) мы будем пользоваться следующими правилами.

- Если область значений какого-либо индекса - это множество, состоящее из одного элемента, мы будем опускать соответствующий индекс.
- Если $n = 1$, то $\sigma_s$ - тождественное преобразование. Это преобразование можно не указывать в выражении.

*определённое равенством*

$$(a_0 \otimes ... \otimes a_n, \sigma) \circ (f_1 \otimes ... \otimes f_n) = a_0 \sigma(f_1) a_1 ... a_{n-1} \sigma(f_n) a_n \qquad (5.6)$$
$$a_0, ..., a_n \in A \quad \sigma \in S_n \quad f_1, ..., f_n \in \mathcal{L}(A_n; A)$$

*является представлением*[6] *алгебры* $A^{\otimes n+1} \times S^n$ *в* $D$*-модуле* $\mathcal{L}(A^n; A)$*.*

*Доказательство.* Согласно теоремам 3.5, 5.6, $n$-линейное отображение можно представить в виде суммы слагаемых (5.1), где $f_i$, $i = 1, ..., n$, генераторы представления (3.3). Запишем слагаемое $s$ выражения (5.5) в виде

$$b_1 \sigma(I_{1\cdot s} \circ x_1) c_1 b_2 ... c_{n-1} b_n \sigma(I_{n\cdot s} \circ x_n) c_n \qquad (5.7)$$

где

$$b_1 = f^n_{s\cdot 0} \quad b_2 = ... = b_n = e \quad c_1 = f^n_{s\cdot 1} \quad ... \quad c_n = f^n_{s\cdot n}$$

Положим в равенстве (5.7)

$$f_i = \sigma^{-1}(b_i) I_{i\cdot s} \sigma^{-1}(c_i) \quad i = 1, ..., n$$

Следовательно, согласно теореме 5.3, отображение (5.6) является преобразованием модуля $\mathcal{L}(A^n; A)$. Для данного тензора $c \in A^{\otimes n+1}$ и заданной перестановки $\sigma \in S_n$, преобразование $h(c, \sigma)$ является линейным преобразованием модуля $\mathcal{L}(A^n; A)$ согласно теореме 5.5. Согласно теореме 5.3, отображение (5.6) является линейным отображением. Согласно определению [3]-2.1.4 отображение (5.6) является представлением алгебры $A^{\otimes n+1} \times S^n$ в модуле $\mathcal{L}(A^n; A)$. □

**Теорема 5.9.** *Рассмотрим* $D$*-алгебру* $A$*. Представление*

$$h : A^{\otimes n+1} \times S_n \to {}^*\mathcal{L}(A^n; A)$$

*алгебры* $A^{\otimes n+1}$ *в модуле* $\mathcal{L}(A^n; A)$*, определённое равенством* (5.6) *позволяет отождествить тензор* $d \in A^{\otimes n+1}$ *и перестановку* $\sigma \in S^n$ *с отображением*

$$(d, \sigma) \circ (f_1, ..., f_n) \quad f_i = \delta \in \mathcal{L}(A; A) \qquad (5.8)$$

*где* $\delta \in \mathcal{L}(A; A)$ *- тождественное отображение.*

*Доказательство.* Если в равенстве (5.2) мы положим $f_i = \delta$, $d = a_0 \otimes ... \otimes a_n$, то равенство (5.2) приобретает вид

$$\begin{aligned} ((a_0 \otimes ... \otimes a_n, \sigma) \circ (\delta, ..., \delta)) \circ (x_1, ..., x_n) &= a_0\ (\delta \circ x_1)\ ...\ (\delta \circ x_n)\ a_n \\ &= a_0\ x_1\ ...\ x_n\ a_n \end{aligned} \qquad (5.9)$$

Если мы положим

$$\begin{aligned} &((a_0 \otimes ... \otimes a_n, \sigma) \circ (\delta, ..., \delta)) \circ (x_1, ..., x_n) \\ &= (a_0 \otimes ... \otimes a_n, \sigma) \circ (\delta \circ x_1, ..., \delta \circ x_n) \\ &= (a_0 \otimes ... \otimes a_n, \sigma) \circ (x_1, ..., x_n) \end{aligned} \qquad (5.10)$$

то сравнение равенств (5.9) и (5.10) даёт основание отождествить действие тензора $d = a_0 \otimes ... \otimes a_n$ и перестановки $\sigma \in S^n$ с отображением (5.8). □

[6] Определение представления $\Omega$-алгебры дано в определении [3]-2.1.4.

Вместо записи $(a_0 \otimes ... a_n, \sigma)$ мы также будем пользоваться записью

$$a_0 \otimes_{\sigma(1)} ... \otimes_{\sigma(n)} a_n$$

если мы хотим явно указать какой аргумент становится на соответствующее место. Например, следующие выражения эквивалентны

$$(a_0 \otimes a_1 \otimes a_2 \otimes a_3, (2,1,3)) \circ (x_1, x_2, x_3) = a_0 x_2 a_1 x_1 a_2 x_3 a_3$$

$$(a_0 \otimes_2 a_1 \otimes_1 a_2 \otimes_3 a_3) \circ (x_1, x_2, x_3) = a_0 x_2 a_1 x_1 a_2 x_3 a_3$$

## 6. Полилинейное отображение в свободную конечно мерную ассоциативную алгебру

**Теорема 6.1.** *Пусть $A$ - свободная конечно мерная ассоциативная алгебра над кольцом $D$. Пусть $\overline{\overline{I}}$ - базис алгебры $\mathcal{L}(A;A)$. Пусть $\overline{\overline{e}}$ - базис алгебры $A$ над кольцом $D$.* **Стандартное представление полилинейного отображения в ассоциативную алгебру** *имеет вид*

$$f \circ (a_1, ..., a_n) = f^{i_0 ... i_n}_{t \cdot k_1 ... k_n} \ \overline{e}_{i_0} \ \sigma_t(I_{k_1} \circ a_1) \ \overline{e}_{i_1} \ ... \ \sigma_t(I_{k_n} \circ a_n) \ \overline{e}_{i_n} \tag{6.1}$$

*Индекс $t$ нумерует всевозможные перестановки $\sigma_t$ множества переменных $\{a_1, ..., a_n\}$. Выражение $f^{i_0 ... i_n}_{t \cdot k_1 ... k_n}$ в равенстве (6.1) называется* **стандартной компонентой полилинейного отображения** *$f$.*

*Доказательство.* Мы изменим индекс $s$ в равенстве (5.5) таким образом, чтобы сгруппировать слагаемые с одинаковым набором генераторов $I_k$. Выражение (5.5) примет вид

$$a = f^n_{k_1 ... k_n \cdot s \cdot 0} \ \sigma_s(I_{k_1 \cdot s} \circ a_1) \ f^n_{k_1 ... k_n \cdot s \cdot 1} \ ... \ \sigma_s(I_{k_n \cdot s} \circ a_n) \ f^n_{k_1 ... k_n \cdot s \cdot n} \tag{6.2}$$

Мы предполагаем, что индекс $s$ принимает значения, зависящие от $k_1$, ..., $k_n$. Компоненты полилинейного отображения $f$ имеют разложение

$$f^n_{k_1 ... k_n \cdot s \cdot p} = \overline{e}_i f^{ni}_{k_1 ... k_n \cdot s \cdot p} \tag{6.3}$$

относительно базиса $\overline{\overline{e}}$. Если мы подставим (6.3) в (5.5), мы получим

$$a = f^{nj_1}_{k_1 ... k_n \cdot s \cdot 0} \ \overline{e}_{j_1} \ \sigma_s(I_{k_1} \circ a_1) \ f^{nj_2}_{k_1 ... k_n \cdot s \cdot 1} \ \overline{e}_{j_2} \ ... \ \sigma_s(I_{k_n} \circ a_n) \ f^{nj_n}_{k_1 ... k_n \cdot s \cdot n} \ \overline{e}_{j_n} \tag{6.4}$$

Рассмотрим выражение

$$f^{j_0 ... j_n}_{t \cdot k_1 ... k_n} = f^{nj_1}_{k_1 ... k_n \cdot s \cdot 0} \ ... f^{nj_n}_{k_1 ... k_n \cdot s \cdot n} \tag{6.5}$$

В правой части подразумевается сумма тех слагаемых с индексом $s$, для которых перестановка $\sigma_s$ совпадает. Каждая такая сумма будет иметь уникальный индекс $t$. Подставив в равенство (6.4) выражение (6.5) мы получим равенство (6.1). □

**Теорема 6.2.** *Пусть $A$ - свободная конечно мерная ассоциативная алгебра над кольцом $D$. Пусть $\overline{\overline{e}}$ - базис алгебры $A$ над кольцом $D$. Полилинейное отображение (5.4) можно представить в виде $D$-значной формы степени $n$ над кольцом $D$*[7]

$$f(a_1, ..., a_n) = a_1^{i_1} ... a_n^{i_n} f_{i_1 ... i_n} \tag{6.6}$$

[7] Теорема доказана по аналогии с теоремой в [2], с. 107, 108

*где*

$$a_j = \overline{e}_i a^i_j \tag{6.7}$$
$$f_{i_1...i_n} = f \circ (\overline{e}_{i_1}, ..., \overline{e}_{i_n})$$

*и величины $f_{i_1...i_n}$ являются координатами $D$-значного ковариантнго тензора.*

*Доказательство.* Согласно определению 1.4 равенство (6.6) следует из цепочки равенств

$$f \circ (a_1, ..., a_n) = f \circ (\overline{e}_{i_1} a_1^{i_1}, ..., \overline{e}_{i_n} a_n^{i_n}) = a_1^{i_1} ... a_n^{i_n} f \circ (\overline{e}_{i_1}, ..., \overline{e}_{i_n})$$

Пусть $\overline{\overline{e}}'$ - другой базис. Пусть

$$\overline{e}'_i = \overline{e}_j h^j_i \tag{6.8}$$

преобразование, отображающее базис $\overline{\overline{e}}$ в базис $\overline{\overline{e}}'$. Из равенств (6.8) и (6.7) следует

$$\begin{aligned} f'_{i_1...i_n} &= f \circ (\overline{e}'_{i_1}, ..., \overline{e}'_{i_n}) \\ &= f \circ (\overline{e}_{j_1} h^{j_1}_{i_1}, ..., \overline{e}'_{j_n} h^{j_n}_{i_n}) \\ &= h^{j_1}_{i_1} ... h^{j_n}_{i_n} f \circ (\overline{e}_{j_1}, ..., \overline{e}_{j_n}) \\ &= h^{j_1}_{i_1} ... h^{j_n}_{i_n} f_{j_1...j_n} \end{aligned} \tag{6.9}$$

Из равенства (6.9) следует тензорный закон преобразования координат полилинейного отображения. Из равенства (6.9) и теоремы [4]-2.1.4 следует, что значение отображения $f \circ (\overline{a}_1, ..., \overline{a}_n)$ не зависит от выбора базиса. □

Полилинейное отображение (5.4) **симметрично**, если

$$f \circ (a_1, ..., a_n) = f \circ (\sigma(a_1), ..., \sigma(a_n))$$

для любой перестановки $\sigma$ множества $\{a_1, ..., a_n\}$.

**Теорема 6.3.** *Если полилинейное отображение $f$ симметрично, то*

$$f_{i_1,...,i_n} = f_{\sigma(i_1),...,\sigma(i_n)} \tag{6.10}$$

*Доказательство.* Равенство (6.10) следует из равенства

$$\begin{aligned} a_1^{i_1} ... a_n^{i_n} f_{i_1...i_n} &= f \circ (a_1, ..., a_n) \\ &= f \circ (\sigma(a_1), ..., \sigma(a_n)) \\ &= a_1^{i_1} ... a_n^{i_n} f_{\sigma(i_1)...\sigma(i_n)} \end{aligned}$$

□

Полилинейное отображение (5.4) **косо симметрично**, если

$$f \circ (a_1, ..., a_n) = |\sigma| f \circ (\sigma(a_1), ..., \sigma(a_n))$$

для любой перестановки $\sigma$ множества $\{a_1, ..., a_n\}$. Здесь

$$|\sigma| = \begin{cases} 1 & \text{перестановка } \sigma \text{ чётная} \\ -1 & \text{перестановка } \sigma \text{ нечётная} \end{cases}$$

**Теорема 6.4.** *Если полилинейное отображение $f$ косо симметрично, то*

$$f_{i_1,...,i_n} = |\sigma| f_{\sigma(i_1),...,\sigma(i_n)} \tag{6.11}$$

*Доказательство.* Равенство (6.11) следует из равенства

$$\begin{aligned} a_1^{i_1}...a_n^{i_n}f_{i_1...i_n} =&f\circ(a_1,...,a_n)\\ =&|\sigma|f\circ(\sigma(a_1),...,\sigma(a_n))\\ =&a_1^{i_1}...a_n^{i_n}|\sigma|f_{\sigma(i_1)...\sigma(i_n)} \end{aligned}$$

□

**Теорема 6.5.** *Пусть $A$ - свободная конечно мерная ассоциативная алгебра над кольцом $D$. Пусть $\overline{\overline{I}}$ - базис алгебры $\mathcal{L}(A;A)$. Пусть $\overline{\overline{e}}$ - базис алгебры $A$ над кольцом $D$. Пусть полилинейное над кольцом $D$ отображение $f$ порождено набором отображений $(I_{k_1},...,I_{k_n})$. Координаты отображения $f$ и его компоненты относительно базиса $\overline{\overline{e}}$ удовлетворяют равенству*

$$f_{l_1...l_n}=f^{i_0...i_n}_{t\cdot k_1...k_n}I_{k_1\cdot}{}^{j_1}_{l_1}...I_{k_n\cdot}{}^{j_n}_{l_n}C^{k_1}_{i_0\sigma_t(j_1)}C^{l_1}_{k_1i_1}...B^{k_n}_{l_{n-1}\sigma_t(j_n)}C^{l_n}_{k_ni_n}\overline{e}_{l_n} \tag{6.12}$$

$$f^p_{l_1...l_n}=f^{i_0...i_n}_{t\cdot k_1...k_n}I_{k_1\cdot}{}^{j_1}_{l_1}...I_{k_n\cdot}{}^{j_n}_{l_n}C^{k_1}_{i_0\sigma_t(j_1)}C^{l_1}_{k_1i_1}...C^{k_n}_{l_{n-1}\sigma_t(j_n)}C^{p}_{k_ni_n} \tag{6.13}$$

*Доказательство.* В равенстве (6.1) положим

$$I_{k_i}\circ a_i=\overline{e}_{j_i}I_{k_i\cdot}{}^{j_i}_{l_i}a_i^{l_i}$$

Тогда равенство (6.1) примет вид

$$\begin{aligned} f\circ(a_1,...,a_n)=&f^{i_0...i_n}_{t\cdot k_1...k_n}\ \overline{e}_{i_0}\sigma_t(a_1^{l_1}I_{k_1\cdot}{}^{j_1}_{l_1}\overline{e}_{j_1})\overline{e}_{i_1}...\sigma_t(a_n^{l_n}I_{k_n\cdot}{}^{j_n}_{l_n}\overline{e}_{j_n})\overline{e}_{i_n}\\ =&a_1^{l_1}...a_n^{l_n}f^{i_0...i_n}_{t\cdot k_1...k_n}I_{k_1\cdot}{}^{j_1}_{l_1}...I_{k_n\cdot}{}^{j_n}_{l_n}\overline{e}_{i_0}\sigma_t(\overline{e}_{j_1})\overline{e}_{i_1}...\sigma_t(\overline{e}_{j_n})\overline{e}_{i_n}\\ =&a_1^{l_1}...a_n^{l_n}f^{i_0...i_n}_{t\cdot k_1...k_n}I_{k_1\cdot}{}^{j_1}_{l_1}...I_{k_n\cdot}{}^{j_n}_{l_n}C^{k_1}_{i_0\sigma_t(j_1)}C^{l_1}_{k_1i_1}\\ &...C^{k_n}_{l_{n-1}\sigma_t(j_n)}C^{l_n}_{k_ni_n}\overline{e}_{l_n} \end{aligned} \tag{6.14}$$

Из равенства (6.6) следует

$$f\circ(a_1,...,a_n)=\overline{e}_pf^p_{i_1...i_n}a_1^{i_1}...a_n^{i_n} \tag{6.15}$$

Равенство (6.12) следует из сравнения равенств (6.14) и (6.6). Равенство (6.13) следует из сравнения равенств (6.14) и (6.15). □

## 7. Список литературы

# 8. Предметный указатель

# 9. Специальные символы и обозначения